\documentclass[onecolumn,reqno]{amsart}

\usepackage{amsmath,amssymb,amsthm,times,amsfonts}
\usepackage {graphicx, color}

\newtheorem{theorem}{Theorem}[section]
\newtheorem{lemma}[theorem]{Lemma}
\newtheorem{remark}[theorem]{Remark}
\newtheorem{corollary}[theorem]{Corollary}
\newtheorem {proposition}[theorem]{Proposition}
\newtheorem {definition}[theorem]{Definition}

\usepackage[pagewise]{lineno}
\def \RR{\mathbb R}

\def\({\left(}  \def\){\right)}
\def\[{\left[}  \def\]{\right]}

\def \beq {\begin {equation}}
\def \eeq {\end{equation}}
\def \OL {\overline}

\def \W {\widetilde}

\begin {document}

	\title [Universal weighted bounds]{ Bounds on energy and potentials of discrete measures on the sphere}

\author[S. Borodachov]{S. V. Borodachov}
\address{Department of Mathematics, Towson University, 7800 York Rd, Towson, MD, 21252, USA}
\email{sborodachov@towson.edu}

\author[P. Boyvalenkov]{P. G. Boyvalenkov}
\address{ Institute of Mathematics and Informatics, Bulgarian Academy of Sciences,
8 G Bonchev Str., 1113  Sofia, Bulgaria}
\email{peter@math.bas.bg}

\author[P. Dragnev]{P. D. Dragnev}
\address{ Department of Mathematical Sciences, Purdue University \\
Fort Wayne, IN 46805, USA }
\email{dragnevp@pfw.edu}

\author[D. Hardin]{D. P. Hardin}
\address{ Center for Constructive Approximation, Department of Mathematics \\
Vanderbilt University, Nashville, TN 37240, USA }
\email{doug.hardin@vanderbilt.edu}

\author[E. Saff]{E. B. Saff}
\address{ Center for Constructive Approximation, Department of Mathematics \\
Vanderbilt University, Nashville, TN 37240, USA }
\email{edward.b.saff@vanderbilt.edu}

\author[M. Stoyanova]{M. M. Stoyanova} 
\address{ Faculty of Mathematics and Informatics, Sofia University ``St. Kliment Ohridski"\\
5 James Bourchier Blvd., 1164 Sofia, Bulgaria}
\email{stoyanova@fmi.uni-sofia.bg}
	
\maketitle


 \begin{abstract}  
We establish upper and lower universal bounds for potentials of weighted designs on the sphere $\mathbb{S}^{n-1}$ that depend only on quadrature nodes and weights derived from the design structure. Our bounds hold for a large class of potentials that includes absolutely monotone functions. The classes of spherical designs attaining these bounds are characterized. Additionally, we study the problem of constrained energy minimization for Borel probability measures on $\mathbb{S}^{n-1}$ and apply it to optimal distribution of charge supported at a given number of points on the sphere. In particular, our results apply to $p$-frame energy.
\end {abstract}

\section {Introduction and review of known results}

Universal upper and lower bounds for potentials of (equi-weighted) spherical designs with a characterization of extremal spherical designs as well as a closely related question of describing all universal maxima (resp. minima) of the extremal spherical designs have recently received a noticeable attention in mathematical literature \cite {Sto1975circle,Sto1975,HarKenSaf2013,NR1,NR2,B2,B,Bor-new,Bor-2,Bor-FL,BDHSS-JMAA,BDHSS-Sharp}. The first part of our paper is devoted to extending these results to the case of weighted spherical designs and their associated potentials (see Theorems \ref{stiff} and \ref{weakly2msharp}, two of the main results of the paper).  

The problem of minimizing the (equi-weighted) energy of a point configuration on the sphere of a fixed cardinality was studied by a large number of authors, see \cite {Y,A2,CohKum2007,BDHSS,BDHSS-AMP,BDHSS-DCC19,BDHSSMathComp,BilGlaPar2022} and references therein. See also  \cite {BHS} for more extensive reviews on this problem. In particular, Bilyk et al \cite {BilGlaPar2022} minimized the energy of a Borel probability measure with respect to the $p$-frame potential showing the optimality of discrete distributions supported at points of tight designs (for certain ranges of the power $p$). Motivated by applications to $p$-frame energy, in the second part of this paper (Section \ref {7c}, and in particular Theorems \ref{genframe} and \ref{antipodal} and Corollaries \ref{nonantipodal} and \ref{antipodal2N} as the other main results), we minimize the energy of a Borel charge distribution subject to additional constraints. Placing additional constraints on the distribution allowed us to extend the range of potential functions under consideration.

Denote by $\mathbb{S}^{n-1}:=\{(t_1,\ldots,t_n)\in \RR^n : t_1^2+\ldots+t_n^2=1\}$ the unit sphere in $\RR^n$ and let $\mathcal J_{n,N}$, $n\geq 2$, $N\geq 1$, be the set of all ordered pairs $(C,W)$, where $C=\{x_1,\ldots,x_N\}$ is a configuration of pairwise distinct points on the sphere $\mathbb{S}^{n-1}$ (called a spherical code) and $W=(w_1,\ldots,w_N)\in (\RR_{>0})^N$ is a vector such that $w_1+\ldots+w_N=1$ (called a normalized vector of weights). Each weighted code $(C,W)$ can be considered as a mapping $w:C\to\mathbb R_{>0}$ defined by $w(x_i)=w_i$, $i=1,\ldots,N$. We will also use the notation $(C,W)=(x_1,\ldots,x_N;w_1,\ldots,w_N)$ and call $(C,W)$ a {\it weighted code}. 

There is a one-to-one correspondence between weighted codes $(C,W)$ in $\mathcal J_{n,N}$ and discrete probability measures $$
\kappa(C,W):=\sum\limits_{i=1}^{N}w_i\delta_{x_i}, 
$$
where $\delta_{x_i}$ is the Dirac delta-function in $\RR^n$ with mass $1$ at $x_i$. There is also a one-to-one correspondence between the set $\mathcal J_{n,N}$ and the set of cubature formulas for integration over the uniform probability distribution on $\mathbb{S}^{n-1}$ which have $N$ pairwise distinct nodes and strictly positive weights and are exact on constant functions. Any multiset on $\mathbb{S}^{n-1}$ (a spherical code where points with distinct indices may coincide) gives a rise to an element $(C,W)$ of $\mathcal J_{n,N}$ for some $N$ where all components of the vector of weights $W$ are rational.

Any function $f:[-1,1]\to (-\infty,\infty]$ finite and differentiable on $(-1,1)$ such that $\lim\limits_{t\to -1^+}f(t)=f(-1)$ and $\lim\limits_{t\to 1^-}f(t)=f(1)$ will be called an {\it admissible potential function}.
Given a weighted code $(C,W)\in \mathcal J_{n,N}$, we define its (polarization) potential with respect to $x\in \mathbb{S}^{n-1}$ as
$$
U_f(x,C;W):=\sum\limits_{i=1}^{N}w_i f(x\cdot x_i){\color{black} \left( =\int_{\mathbb{S}^{n-1}} f(x\cdot y)\, d\kappa(C,W)(y)\right)},
$$
where $x\cdot y:=t_1s_1+\ldots+t_ns_n$ denotes the standard dot product of vectors $x=(t_1,\ldots,t_n)$ and $y=(s_1,\ldots,s_n)$ in $\RR^n$.

Our work is related to the following problem. Given a weighted code $(C,W)\in \mathcal J_{n,N}$ and an admissible potential function $f$, find 
\begin {equation}\label{inf}
\inf\limits_{x\in \mathbb{S}^{n-1}}U_f(x,C;W)
\end {equation}
and points $x^\ast\in \mathbb{S}^{n-1}$ that attain the infimum of \eqref{inf}. If  $f$ is finite at $t=1$ and $t=-1$, one may consider a similar problem of maximizing the potential.

In the equi-weighted case; i.e., when $W=\(1/N,\ldots,1/N\)$, problem \eqref {inf} was solved for several important classes of spherical designs. Recall that a spherical code $C=\{x_1,\ldots,x_N\}\subset \mathbb{S}^{n-1}$ is called a spherical $m$-design if for every polynomial $p$ in $n$ variables of degree at most $m$,
$$
\frac {1}{N}\sum\limits_{i=1}^{N}p(x_i)=\int\limits_{\mathbb{S}^{n-1}}p(x)\ \! d\sigma _n (x),
$$
where $\sigma_n$ is the area measure on $\mathbb{S}^{n-1}$ normalized to be a probability measure.

Following \cite {CohKum2007}, we call a spherical code $C$ {\it $m$-sharp}, $m\geq 1$, if $C$ is a $(2m-1)$-design with $m$ values of the dot product occurring between distinct points of $C$ and {\it strongly $m$-sharp} if, in addition, $C$ is a $2m$-design. A list of known sharp and strongly sharp configurations on the Euclidean sphere can be found, for example, in \cite [Table 1]{CohKum2007}. A spherical code $C\subset \mathbb{S}^{n-1}$ is called {\it $m$-stiff} if $C$ is a $(2m-1)$-design and there is a point on $\mathbb{S}^{n-1}$ forming $m$ distinct values of the dot product with points of $C$. 

Recall also that a function $f:[-1,1)\to\RR$, $\RR=(-\infty,\infty)$, is called {\it absolutely monotone} if $f^{(k)}\geq 0$ on $[-1,1)$ for all $k\geq 0$ and {\it strictly absolutely monotone} if, in addition, $f^{(k)}>0$ on $(-1,1)$ for all $k\geq 0$. A point $z\in \mathbb{S}^{n-1}$ will be called a universal minimum (maximum) of a spherical code $C\subset \mathbb{S}^{n-1}$ if for every admissible potential function $f$ absolutely monotone on $[-1,1)$, the equi-weighted potential
$$
U_f(x,C):=\frac {1}{N}\sum\limits_{i=1}^{N}f(x\cdot x_i)
$$
attains its absolute minimum (maximum) over $\mathbb{S}^{n-1}$ at $z$.
In the earlier papers \cite {Sto1975circle,Sto1975,NR1,NR2} which considered the Riesz potential and in recent works \cite {B2,BorESI,B,Bor-2} universal minima and maxima of the regular $N$-gon on $\mathbb{S}^1$ and of the regular simplex, cross-polytope, and cube on $\mathbb{S}^{n-1}$ were found\footnote {Absolute maxima in the case of a cube were only found for the Riesz potential, see \cite {NR2}. The proof in \cite {NR2} uses induction on dimension and can be extended to absolutely monotone potential functions (to appear in our upcoming work).}. Some of these four basic codes are sharp antipodal, some are strongly sharp, and some are stiff. For any $m$-stiff code, the set of its universal minima is exactly the set of points on $\mathbb{S}^{n-1}$ forming $m$ distinct dot products with points of the code (cf. \cite {BorESI,Bor-2}). In fact, many known sharp codes listed in \cite [Table 1]{CohKum2007} are also stiff (cf. \cite {BDHSS-Sharp}). At the same time, several interesting stiff codes are not sharp (e.g., cube on $\mathbb{S}^{n-1}$, $n\geq 3$, symmetrized regular simplex on $\mathbb{S}^{n-1}$, $n\geq 5$ odd, and the $24$-cell on $\mathbb{S}^3$ with more examples given in \cite {Bor-2}). Sharp antipodal and strongly sharp codes form the class of tight designs as defined in \cite {DGS}. For any strongly sharp code, the set of its universal minima is the set of antipodes of points of the code and the set of its universal maxima is the code itself (cf. \cite {BorESI,B}). Finally, for any sharp antipodal code, the set of its universal maxima is also the code itself \cite {BorESI,B}. Certain remarkable spherical codes are not sharp (for example, the $600$-cell on $\mathbb{S}^3$). Universal maxima of the $600$-cell were shown to be exactly the points of the $600$-cell itself in \cite {BDHSS-JMAA}. The problem about universal minima of the $600$-cell remains open. Certain  sharp antipodal codes that are of great interest in mathematics are not stiff, for example, regular icosahedron and dodecahedron, minimal vectors of the $E_8$ and Leech lattices, and the $2_{41}$ polytope in $\RR^8$. Finding the universal minima for these codes required a more sophisticated technique, see \cite {Bor-new} for the first three codes, \cite {BDHSS-Sharp} for the Leech lattice, and \cite {Bor-FL} for the $2_{41}$ polytope. 

The paper is organized as follows. Weighted spherical designs are introduced in Section 2 and special classes of weighted designs are defined in Section 4. Gauss-type quadratures and associated orthogonal polynomials are reviewed in Section 3. Our universal bounds on potentials of weighted designs are derived in Section 5 and illustrated with examples in Section 6. The constrained (frame) energy minimization results are found in Section 7 with needed auxiliary results presented in an Appendix.

\section {Definition of a weighted spherical design}

Denote by $\mathbb P_n$, $n\geq 0$, the set of all polynomials of one variable of degree at most $n$. 
\begin {definition}\label {const}
{\rm
Let $n\geq 2$, $N\geq 1$ and $m\geq 0$. A weighted code $(C,W)\in \mathcal J_{n,N}$ is called a (weighted) spherical $m$-design if, for any polynomial $p\in \mathbb P_m$,
\begin {equation}\label {polynomial}
\sum\limits_{i=1}^{N}w_ip(x\cdot x_i)=c,\ \ \ x\in \mathbb{S}^{n-1},
\end {equation}
for some constant $c$ independent of $x$.}
\end {definition}

\begin {remark}\label {rc}
{\rm
If \eqref {polynomial} holds for a polynomial $p$, then 
\begin {equation}\label {designequiv}
c=\alpha_0(p):=\int_{-1}^{1}p(t)w_n(t)\ \! dt,
\end {equation}
where $w_n(t)=\gamma_n(1-t^2)^{(n-3)/2}$, $t\in [-1,1]$, and $\gamma_n$ is a positive constant such that $w_n$ is a probability density on $[-1,1]$.
}
\end {remark}
There are several equivalent definitions of a (weighted) spherical design. We will mention one of them. Denote by $P_m^{(n)}$ the $m$-th Gegenbauer polynomial corresponding to the sphere $\mathbb{S}^{n-1}$.
It is orthogonal to the polynomial space $\mathbb P_{m-1}$ with weight $w_n$ and is normalized by setting $P_m^{(n)}(1)=1$. A weighted code $(C,W)\in \mathcal J_{n,N}$ is called a (weighted) spherical $m$-design if
$$
\sum\limits_{i=1}^{N}\sum\limits_{j=1}^{N}w_iw_jP_k^{(n)}(x_i\cdot x_j)=0,\ \ \ 1\leq k\leq m.
$$
The equivalence of this definition to Definition \ref {const} can be proved similarly to the equi-weighted case which can be found in, e.g., \cite [Chapter 5]{BHS}.

Many spherical designs are known when $W=(1/N,\ldots,1/N)$; i.e., ({equi-weighted) spherical designs as defined in \cite {DGS}. Widely known examples are sets of vertices of a regular $N$-gon inscribed in $\mathbb{S}^1$, of a regular cross-polytope, simplex, and cube inscribed in $\mathbb{S}^{n-1}$, of the regular icosahedron and dodecahedron on $\mathbb{S}^2$, of the $24$-cell, $120$-cell, and $600$-cell on $\mathbb{S}^3$, and sets of minimal vectors of $E_8$ and Leech lattices.
As an example of a weighted design, we mention the weighted code $(C,W)\in J_{3,5}$, where $C=C_1 \cup C_2$ consists of $C_1=\{ (0,0,1)\}$ with a weight $1/4$ and $C_2=\{ (\pm 2/3,\pm 2/3,-1/3)\}$, four points with equal weights $3/16$. This is the only weighted code $(C,W)\in \mathcal J_{3,5}$, where $C$ is a square-based pyramid on $\mathbb{S}^2$, which is a weighted spherical $2$-design. As other examples, we mention pairs cube--cross-polytope, icosahedron--dodecahedron, $120$-cell--$600$-cell, and the $56$-point sharp code--roots of $E_7$ lattice, which with certain weights form weighted $5$-, $9$-, $19$-, and $7$-designs, respectively, see \cite {GoeSei1981,Yud2005,BorBoyDraenergy}. More examples are given in \cite {GoeSei1981} and in this paper.

\section {Quadratures of the highest algebraic degree of precision}

Before stating the main definitions and results, we will review the necessary background material.
Let $-1<\alpha_1<\ldots<\alpha_m<1$ be the zeros of the $m$-th Gegenbauer polynomial $P_m^{(n)}$. We will also denote it by $P_m^{0,0}$. Let also $(\rho_1,\ldots,\rho_m)$ be the unique vector of coefficients such that the quadrature
\begin {equation}\label {quad}
\alpha_0(p)=\int_{-1}^{1}p(t)w_n(t)\ \!dt\approx \sum\limits_{i=1}^{m}\rho_ip(\alpha_i)
\end {equation}
is exact for all polynomials $p$ of degree up to $2m-1$. Quadrature \eqref {quad} is known as the Gauss-Gegenbauer quadrature.

Let $-1<\gamma_1<\ldots<\gamma_m<1$ be the zeros of the $m$-th Jacobi polynomial $J_m^{(\frac {n-1}{2},\frac {n-3}{2})}$; i.e. the polynomial of degree $m$ orthogonal to $\mathbb P_{m-1}$ with weight $(1-t)w_n(t)$. We will also denote it by $P_m^{1,0}$.  Let $(\sigma_1,\ldots,\sigma_{m+1})$ be the unique vector of coefficients such that the quadrature
\begin {equation}\label {quad1}
\alpha_0(p)=\int_{-1}^{1}p(t)w_n(t)\ \!dt\approx \sum\limits_{i=1}^{m}\sigma_ip(\gamma_i)+\sigma_{m+1}p(1)
\end {equation}
is exact for all polynomials $p$ of degree up to $2m$. 

The numbers $-1<-\gamma_m<\ldots<-\gamma_1<1$ are the zeros of the $m$-th degree Jacobi polynomial $J_m^{(\frac {n-3}{2},\frac {n-1}{2})}$ (it is orthogonal to $\mathbb P_{m-1}$ with weight $(1+t)w_n(t)$). We will denote it by $P_m^{0,1}$.  The quadrature
\begin {equation}\label {quad2}
\alpha_0(p)=\int_{-1}^{1}p(t)w_n(t)\ \!dt\approx \sigma_{m+1}p(-1)+\sum\limits_{i=1}^{m}\sigma_ip(-\gamma_i)
\end {equation}
is exact for all polynomials $p$ of degree up to $2m$.

Let $-1<\W\alpha_1<\ldots<\W\alpha_m<1$ be the zeros of the $m$-th Gegenbauer polynomial $P_m^{(n+2)}$ corresponding to the sphere $S^{n+1}$ denoted also by $P_m^{1,1}$. It is orthogonal to $\mathbb P_{m-1}$ with weight $w_{n+2}(t)=(1-t^2)w_n(t)$.
Let $(\theta_0,\ldots,\theta_{m+1})$ be the unique vector of coefficients such that the quadrature
\begin {equation}\label {quad3}
\alpha_0(p)=\int_{-1}^{1}p(t)w_n(t)\ \!dt\approx \theta_0p(-1)+\sum\limits_{i=1}^{m}\theta_ip(\W\alpha_i)+\theta_{m+1}p(1)
\end {equation}
is exact for all polynomials $p$ of degree up to $2m+1$. Levenshtein calls $P_m^{\mu,\nu}$, $\mu,\nu\in \{0,1\}$, adjacent polynomials.

Quadratures \eqref {quad}--\eqref {quad3} are special cases of the Gauss-Jacobi mechanical quadrature \textcolor{black}{(see, e.g. \cite[Chapter 3.4]{Sze1975})}. It is convenient to describe them as one quadrature with parameters.
Let $\mu=0$ or $1$ and $\nu=0$ or $1$. The polynomials $P_m^{\mu,\nu}$ are normalized by requiring that $P_m^{\mu,\nu}(1)=1$. 
For each pair $(\mu,\nu)$, the corresponding quadrature \eqref {quad}--\eqref {quad3} has the form
\begin {equation}\label {quadL}
\int_{-1}^{1}p(t)w_n(t)\ \! dt\approx\sum\limits_{k=1-\nu}^{m+\mu}\tau_kp(\lambda_k),
\end {equation}
where the nodes $-1\leq\lambda_{1-\nu}<\ldots<\lambda_{m+\mu}\leq1$ are the zeros of the polynomial $P_m^{\mu,\nu}(t)(1-t)^\mu(1+t)^\nu$ and $\tau_{1-\nu},\ldots,\tau_{m+\mu}$ are the unique weights such that quadrature \eqref {quadL} is exact on $\mathbb P_L$ for $L:=2m-1+\mu+\nu$. The nodes $\lambda_j$ and weights $\tau_j$ are functions of $(\mu,\nu)$; however, we omitted indices $\mu$ and $\nu$ in their notation. Quadrature \eqref {quad} corresponds to the case $\mu=\nu=0$ with $\tau_j=\rho_j$, $j=1,\ldots,m$, quadrature \eqref {quad1} corresponds to the case $\mu=1$, $\nu=0$ with $\tau_j=\sigma_j$, $j=1,\ldots,m+1$, quadrature \eqref {quad2} corresponds to the case $\mu=0$, $\nu=1$ with $\tau_j=\sigma_{m+1-j}$, $j=0,\ldots,m$, and quadrature \eqref {quad3} corresponds to the case $\mu=\nu=1$ with $\tau_j=\theta_j$, $j=0,\ldots,m+1$.

The weights of every quadrature \eqref {quadL} (and, hence, of each quadrature \eqref {quad}--\eqref {quad3}) are strictly positive, see \cite {Sze1975}. Indeed,
denote by $\varphi_j$ the polynomial of degree $m-1$ such that $\varphi_j(\lambda_k)=0$ for $k\in \{1,\ldots,m\}\setminus\{j\}$ and $\varphi_j(\lambda_j)=1$, $j=1,\ldots,m$, and let 
\begin {equation}\label {Rj}
R_j(t):=(\varphi_j(t))^2\frac {(1-t)^\mu(1+t)^\nu}{(1-\lambda_j)^\mu(1+\lambda_j)^{\nu}}.
\end {equation}
Denote $\varphi(t)=(t-\lambda_1)\cdots (t-\lambda_m)$ and let
\begin {equation}\label {R0}
R_{m+1}(t):=\(\frac {\varphi(t)}{\varphi(1)}\)^2\(\frac {1+t}{2}\)^\nu\ \ \ \text {if $\mu=1$ and}\ \ \ R_0(t):=\(\frac {\varphi(t)}{\varphi(-1)}\)^2\(\frac {1-t}{2}\)^\mu
\end {equation}
if $\nu=1$. The polynomials $R_j$
are positive on $[-1,1]$ except for finitely many points and have degree at most $L$. Then quadrature \eqref {quadL} is exact on each $R_j$ and we have
\begin {equation}\label {qe}
\begin {split}
\tau_j&=\sum\limits_{k=1-\nu}^{m+\mu}\tau_kR_j(\lambda_k)
=\int_{-1}^{1}R_j(t)w_n(t)\ \! dt>0,\ \ \ j=1-\nu,\ldots,m+\mu.
\end {split}
\end {equation}


\section {Weakly sharp weighted designs}

We now define the classes of weighted designs, for which the corresponding universal bounds will be sharp. Furthermore, if the derivative of the potential function of the corresponding order is strictly positive, sharpness of the universal bound will imply that the weighted design must be in one of these three classes.

\begin {definition}\label {Dstiff}
{\rm
A weighted code $(C,W)\in \mathcal J_{n,N}$ is called $m$-stiff, $m\geq 1$, if the weighted code $(C,W)$ is a weighted spherical $(2m-1)$-design and there exists a point $z$ on $\mathbb{S}^{n-1}$ that forms at most $m$ distinct values of dot product with all points of the spherical code~$C$.
}
\end {definition}

If it does not cause any confusion, we will omit the word ``weighted" when discussing designs $(C,W)$.

\begin {definition}\label {D2m}
{\rm
A weighted code $(C,W)\in \mathcal J_{n,N}$ is called a {\it weakly sharp $2m$-design}, $m\geq 1$, if in addition to being a $2m$-design, it contains a point $z$ that forms at most $m$ distinct dot products with other points of $C$.
}
\end {definition}

\begin {definition}\label {Dsharp}
{\rm
A weighted code $(C,W)\in \mathcal J_{n,N}$ is called a {\it weakly sharp $(2m+1)$-design}, $m\geq 0$, if it is a $(2m+1)$-design containing a point $z$ that forms at most $m+1$ distinct dot products with other points of $C$ one of which is $-1$.
}
\end {definition}

In the equi-weighted case, Definition \ref {D2m} is equivalent to the definition of a strongly $m$-sharp code and Definition \ref {Dsharp} is equivalent to that of an $(m+1)$-sharp antipodal code. For example, let $C \subset \mathbb{S}^{n-1}$ be an equi-weighted spherical $2m$-design and $z \in C$ be such that $\{z\cdot x : x \in C \setminus \{z\} \}=\{t_1,\ldots,t_m\}$. Using Definition \ref{const} with $p(t)=(t-t_1)^2\ldots(t-t_m)^2$ and $z$, we conclude
that $p_0|{C}|=p(1)$. Now Definition \ref{const} with the same polynomial and any point $y \in C$ implies that $\{y\cdot x : x \in C \setminus \{y\} \}=\{t_1,\ldots,t_m\}$ (the argument for Definition \ref {Dsharp} is similar).
%
%

One can show that in the weighted case, Definition \ref {D2m} does not imply that $C$ is an $m$-distance set, and Definition \ref {Dsharp} does not imply that $C$ is an $(m+1)$-distance set.

We will call point $z$ in Definitions \ref {Dstiff}--\ref {Dsharp} an {\it extremal point} of the stiff or weakly sharp design $C$.
For such a code $C$ define the set
$$
A_{\mu,\nu}:=\begin {cases}
\mathbb{S}^{n-1}, & (\mu,\nu)=(0,0), \cr
C, & (\mu,\nu)=(1,0), \cr
-C, & (\mu,\nu)=(0,1), \cr
C\cap (-C), & (\mu,\nu)=(1,1).\cr
\end {cases}
$$
\begin {remark}
{\rm
Let a weighted code $(C,W)\in \mathcal J_{n,N}$ be a $(2m-1+\mu+\nu)$-design such that there is a point in the set $A_{\mu,\nu}$ forming at most $m+\mu+\nu$ distinct dot products with points of $C$. If $(\mu,\nu)=(0,0)$, the weighted code $(C,W)$ is $m$-stiff, if $(\mu,\nu)=(1,0)$ or $(0,1)$, the weighted code $(C,W)$ is a weakly sharp $2m$-design, and if $(\mu,\nu)=(1,1)$, then $(C,W)$ is a weakly sharp $(2m+1)$-design.
}
\end {remark}

A point on the sphere cannot form ``too few" dot products with points of a given design. The next lemma is a weighted analog to Theorems 6.6 and 6.7 in \cite{DGS} for the equi-weighted case.

\begin {lemma}\label {allofthem1}
Let $n\geq 2$, $m\geq 0$, $N\geq 1$, $\alpha,\beta\in \{0,1\}$, a weighted code $(C,W)\in \mathcal J_{n,N}$ be a $(2m+\alpha+\beta)$-design, and $z\in A_{\alpha,\beta}$ be arbitrary point. Then $z$ must form at least $m+1+\alpha+\beta$ distinct dot products with points of $C$. 
\end {lemma}
\begin {proof}
When $m=0$, we will assume that $\alpha=1$ or $\beta=1$ (otherwise, the assertion is trivial).
Assume that $C$ is a $(2m+\alpha+\beta)$-design and that $z\in A_{\alpha,\beta}$ is an arbitrary point. Assume that $z$ forms a set of $k+\alpha+\beta$ distinct dot products with points $x_1,\ldots,x_N$ from $C$. If $\alpha=1$, we exclude $1$ from this set, if $\beta=1$, we exclude $-1$ from this set, and denote the remaining dot products by $-1\leq t_1<\ldots<t_k\leq 1$. 

Assume to the contrary that $k\leq m$. Let $p(t):=(t+1)^\beta (t-t_1)^2\cdots(t-t_k)^2(1-t)^\alpha$, which is strictly positive except for finitely many points.
Since $(C,W)$ is a $(2m+\alpha+\beta)$-design and ${\rm deg}\ \! p=2k+\alpha+\beta\leq 2m+\alpha+\beta$, in view of \eqref {polynomial} and \eqref {designequiv}, we have
$$
0=\sum\limits_{i=1}^{N}w_i p(z\cdot x_i)=\alpha_0(p)=\int_{-1}^{1}p(t)w_n(t)\ \! dt>0.
$$
This contradiction shows that $k\geq m+1$ proving the lemma.
\end {proof}

Extremal points of a weighted design have the following distributions of dot products with points of the design. 

\begin {lemma}\label {allofthem}
Let $n\geq 2$, $m,N\geq 1$, $\mu,\nu\in \{0,1\}$. If $(C,W)$ is an $L$-design, $L=2m-1+\mu+\nu$, and a point $z\in A_{\mu,\nu}$ forms exactly $m+\mu+\nu$ distinct dot products with points of $C$, then they are $\lambda_{1-\nu},\ldots,\lambda_{m+\mu}$ as defined in \eqref{quadL}, and for each $j=1-\nu,\ldots,m+\mu$, the sum of the weights of points in $C$ forming dot product $\lambda_j$ with $z$ equals~$\tau_j$.
\end {lemma}

\begin {proof}
Assume that $(C,W)$ is an $L$-design and some point $z\in A_{\mu,\nu}$ forms exactly $M:=m+\mu+\nu$ distinct dot products with points of the code $C$. Denote all the dot products by $t_{1-\nu}<\ldots<t_{m+\mu}$. Then $t_0:=-1$ and $t_{m+1}=1$, whenever present. Denote $q(t):=(t-t_1)\cdots(t-t_m)$.  Let $\beta_j$, $j=1-\nu,\ldots,m+\mu$, be the sum of the weights $w_i$ of all points in $C$ forming dot product $t_j$ with $z$. Let $r(t)$ be any polynomial of degree at most $m-1$. Since $(C,W)$ is a $L$-design and the polynomial $p(t)r(t)$, where $p(t):=(1+t)^\nu (t-t_1)\cdots (t-t_m)(1-t)^\mu$, has degree at most $L$, we have
\begin{eqnarray*}
\int_{-1}^{1} q(t)r(t)(1-t)^\mu (1+t)^\nu w_n(t)\ \! dt &=& \int_{-1}^{1}p(t)r(t)w_n(t)\ \! dt =\ \! \alpha_0(pr)\\
&=& \sum\limits_{i=1}^{N}w_i p(z\cdot x_i)r(z\cdot x_i)\ =\ \sum\limits_{j=1-\nu}^{m+\mu}\beta_jp(t_j)r(t_j)\ = \ 0.
\end{eqnarray*}
Then $q(t)$ is the $m$-th degree orthogonal polynomial with the weight $(1-t)^\mu(1+t)^\nu w_n(t)$; i.e., it equals $P_m^{\mu,\nu}$ modulo a constant factor. Recall that all zeros of $P_m^{\mu,\nu}$ lie in $(-1,1)$ \cite {Sze1975}. Hence, $\{t_1,\ldots,t_m\}$ is the set $\{\lambda_1,\ldots,\lambda_m\}$ of zeros of $P_m^{\mu,\nu}$ and $t_j=\lambda_j$, $j=1-\nu,\ldots,m+\mu$.

Since $\lambda_0:=-1$, $\lambda_{m+1}=1$, $(C,W)$ is a $L$-design, and ${\rm deg}\ \! R_j\leq L$, where the polynomials $R_j$, $j=1-\nu,\ldots,m+\mu$, are defined in \eqref {Rj} and \eqref {R0},
by \eqref {qe} we have
\begin {equation*}
\begin {split}
\tau_j&=\int_{-1}^{1}R_j(t)w_n(t)\ \! dt=\alpha_0(R_j)=\sum\limits_{i=1}^{N}w_i R_j(z\cdot x_i)=\beta_jR_j(\lambda_j)=\beta_j;
\end {split}
\end {equation*}
that is, the sum of the weights of points in $C$ that form dot product $\lambda_j$ with $z$ is $\tau_j$, $j=1-\nu,\ldots,m+\mu$.
\end {proof}

\section {Universal bounds on potentials of weighted designs}

In this section we state and prove two of our main results. We start with the following auxiliary statement.

\begin {lemma}\label {bounded}
If $f:(-1,1)\to \RR$ has a bounded below derivative $f''$ on $(-1,1)$, then $f'$ is bounded above on $(-1,0)$ and bounded below on $(0,1)$ while $f$ is bounded below on $(-1,1)$.
\end {lemma}
\begin {proof}
Let $b$ be a lower bound for $f''$ on $(-1,1)$. Then
for any $t\in (-1,0)$, there is $\xi\in (t,0)$ such that
$$
f'(t)-f'(0)=f''(\xi)t\leq bt\leq \max\{-b,0\}
$$
while for any $t\in (0,1)$, there is $\eta\in (0,t)$ such that
$$
f'(t)-f'(0)=f''(\eta)t\geq bt\geq \min\{b,0\}.
$$
Then $f'$ is bounded above on $(-1,0)$, say by $b_1$, and bounded below on $(0,1)$, say by $b_2$.
For any $t\in (-1,0)$, there is $\zeta\in (t,0)$ such that
$$
f(t)-f(0)=f'(\zeta)t\geq b_1t\geq \min\{-b_1,0\}
$$
and for any $t\in (0,1)$, there is $\tau\in (0,t)$ such that
$$
f(t)-f(0)=f'(\tau)t\geq b_2t \geq \min\{b_2,0\}.
$$
Then $f$ is bounded below on $(-1,1)$.
\end {proof}

The following universal lower bound holds for the potential of a given weighted design. We also characterize weighted designs of a given strength that attain this universal lower bound for that strength. Here, $\tau_j$ and $\lambda_j$ are weights and nodes of quadrature \eqref {quadL} for the corresponding pair $(\mu,\nu)$.

\begin {theorem}\label {stiff}
Let $\mu=0$, $\nu=0$ or $1$, $n\geq 2$, $m,N\geq 1$, $f:[-1,1]\to(-\infty,\infty]$ be an admissible potential function such that $f^{(2m+\nu)}\geq 0$ on $(-1,1)$, and a weighted code $(C,W)\in\mathcal J_{n,N}$, be a $(2m-1+\nu)$-design. 

Then 
\begin {equation}\label {potential-lower}
U_f(x,C;W)=\sum\limits_{i=1}^{N}w_if(x\cdot x_i)\geq \sum\limits_{j=1-\nu}^{m}\tau_jf(\lambda_j),\ \ \ x\in \mathbb{S}^{n-1},
\end {equation}
where $\lambda_j$ and $\tau_j$ are the nodes and weights, respectively, of quadrature \eqref {quadL} with $\mu=0$.

If $\nu=0$ and $(C,W)$ is $m$-stiff or $\nu=1$ and $(C,W)$ is a weakly sharp $2m$-design, then equality holds in \eqref {potential-lower} at every point $x\in A_{0,\nu}$ forming $m+\nu$ distinct dot products with points of the code $C$.

If $f^{(2m+\nu)}>0$ on $(-1,1)$ and equality holds in \eqref {potential-lower} at some point $x\in \mathbb{S}^{n-1}$, then (i) the point $x$ lies in $A_{0,\nu}$ and forms exactly $m+\nu$ distinct dot products $\lambda_{1-\nu},\ldots,\lambda_m$ with points of the code $C$, (ii) the weighted code $(C,W)$ is $m$-stiff if $\nu=0$ or is a weakly sharp $2m$-design if $\nu=1$, and (iii) for every $j=1-\nu,\ldots,m$, the sum of the weights of points from $C$ forming dot product $\lambda_j$ with $x$ is exactly~$\tau_j$.
\end {theorem}
We remark that by Lemma \ref {bounded}, we have $f(-1)<\infty$ in the above theorem when $\nu=1$. Indeed, from $f^{(2m+1)} \geq 0$ we conclude inductively that $f' $ is bounded below, which implies that $f$ is bounded above on $(-1,0)$.

\begin {proof}[Proof of Theorem \ref {stiff}]
Let $q\in \mathbb P_{2m-1+\nu}$ be the polynomial interpolating $f$ and $f'$ at points $\lambda_1,\ldots,\lambda_m$ which also interpolates $f$ at $\lambda_0=-1$ if $\nu=1$. For every $t\in (-1,1)$, there is a point $\xi\in (-1,1)$ such that 
\begin {equation}\label {c_1'}
f(t)-q(t)=\frac {f^{(2m+\nu)}(\xi)}{(2m+\nu)!}(t-\lambda_1)^2\cdots (t-\lambda_m)^2(t+1)^\nu.
\end {equation}
If $f(1)<\infty$, then \eqref {c_1'} holds for $t=1$. If $f(-1)<\infty$, then \eqref {c_1'} also holds for $t=-1$.
Since, by assumption, $f^{(2m+\nu)}\geq 0$ on $(-1,1)$, we have $f(t)-q(t)\geq 0$ or $q(t)\leq f(t)$, $t\in [-1,1]$ (these inequalities hold trivially at $t=1$ or $-1$ if $f(1)=\infty$ or $f(-1)=\infty$, respectively). Since $(C,W)$ is a $(2m-1+\nu)$-design and quadrature \eqref {quadL} is exact on $q(t)$, for every point $x\in \mathbb{S}^{n-1}$, using \eqref {designequiv} we have
\begin {equation}\label {pq'}
\begin {split}
U_f(x,C;W)&=\sum\limits_{i=1}^{N}w_if(x\cdot x_i)\geq \sum\limits_{i=1}^{N}w_iq(x\cdot x_i)=\alpha_0(q)\\
&=\sum\limits_{j=1-\nu}^{m}\tau_jq(\lambda_j)=\sum\limits_{j=1-\nu}^{m}\tau_jf(\lambda_j).
\end {split}
\end {equation}

Let now $(C,W)$ be $m$-stiff if $\nu=0$ or a weakly sharp $2m$-design if $\nu=1$. Let also $z\in A_{0,\nu}$ be any point that forms $m+\nu$ distinct values of the dot product with points from $C$ (such $z$ exists by definition). By Lemma \ref {allofthem}, the values of these $m+\nu$ dot products are the zeros $\lambda_{1-\nu}<\ldots<\lambda_m$ of $(1+t)^\nu P_m^{0,\nu}(t)$ and for each $j=1-\nu,\ldots,m$, the sum of the weights of points in $C$ that form dot product $\lambda_j$ with $z$ is $\tau_j$. Then
$$
U_f(z,C;W)=\sum\limits_{i=1}^{N}w_if(z\cdot x_i)=
\sum\limits_{j=1-\nu}^{m}\tau_jf(\lambda_j).
$$

Assume now that $f^{(2m+\nu)}>0$ on $(-1,1)$. Then in view of \eqref {c_1'} we have $f(t)>q(t)$ for $t\in [-1,1]\setminus \{\lambda_{1-\nu},\ldots,\lambda_m\}$. Let $x\in \mathbb{S}^{n-1}$ be any point such that equality holds in \eqref {potential-lower}. Then we have equality in the first line of \eqref {pq'}. Since $f\geq q$ on $[-1,1]$, we cannot have $f(x\cdot  x_i)>q(x\cdot x_i)$ for any $i$ without ruining this equality. Therefore, each dot product $ x\cdot x_i$ is a point of equality of $f$ and $q$; that is, one of the numbers $\lambda_{1-\nu},\ldots,\lambda_m$.

If $\nu=0$, since $(C,W)$ is a $(2m-2)$-design, by Lemma \ref {allofthem1} with $\alpha=\beta=0$ (note that then $A_{\alpha,\beta}=\mathbb{S}^{n-1}$) and $m-1$ instead of $m$, dot product $x\cdot x_i$ must assume at least $m$ distinct values. Then it assumes each of the values $\lambda_{1},\ldots,\lambda_m$.

If $\nu=1$, since $(C,W)$ is a $2m$-design, by Lemma \ref {allofthem1} with $\alpha=\beta=0$, any point on $\mathbb{S}^{n-1}$ must form at least $m+1$ dot products with points from $C$. Therefore, $x$ forms all dot products $\lambda_0,\ldots,\lambda_m$. In particular, $x\in -C=A_{0,1}$.

Thus, in both cases, $x\in A_{0,\nu}$ and dot product $x\cdot x_i$ assumes exactly $m+\nu$ values $\lambda_{1-\nu},\ldots,\lambda_m$. By definition, the $(2m-1+\nu)$-design $(C,W)$ is $m$-stiff if $\nu=0$ or is a weakly sharp $2m$-design if $\nu=1$. By Lemma \ref {allofthem}, the sum of the weights of points in $C$ forming dot product $\lambda_j$ with $x$ is $\tau_j$, $j=1-\nu,\ldots,m$. 
\end {proof}

A similar universal upper bound holds for the potential of a given weighted design. We also characterize weighted designs of a given strength that attain this universal upper bound. Note that $A_{1,\nu}\subset C$, $\nu=0,1$.
\begin {theorem}\label {weakly2msharp}
Let $\mu=1$, $\nu=0$ or $1$, $n\geq 2$, $m,N\geq 1$, $f:[-1,1]\to\RR$ be an admissible potential function such that $f^{(2m+1+\nu)}\geq 0$ on $(-1,1)$, and a weighted code $(C,W)\in\mathcal J_{n,N}$ be a $(2m+\nu)$-design. 

Then 
\begin {equation}\label {potential-upper}
U_f(x,C;W)=\sum\limits_{i=1}^{N}w_if(x\cdot x_i)\leq  \sum\limits_{j=1-\nu}^{m+1}\tau_jf(\lambda_j),\ \ \ x\in \mathbb{S}^{n-1},
\end {equation}
where $\lambda_j$ and $\tau_j$ are the nodes and weights, respectively, of quadrature \eqref {quadL} with $\mu=1$.

If $(C,W)$ is a weakly sharp $(2m+\nu)$-design, then equality holds in \eqref {potential-upper} at every point $x\in A_{1,\nu}$ forming $m+\nu$ distinct dot products with other points of $C$.

If $f^{(2m+1+\nu)}>0$ on $(-1,1)$ and equality holds in \eqref {potential-upper} at some point $x\in \mathbb{S}^{n-1}$, then (i) the point $x$ lies in $A_{1,\nu}$ and forms exactly $m+1+\nu$ distinct dot products, $\lambda_{1-\nu},\ldots,\lambda_{m+1}$, with points of the code $C$, (ii) the weighted code $(C,W)$ is a weakly sharp $(2m+\nu)$-design,  and (iii) for every $j=1-\nu,\ldots,m+1$, the sum of the weights of points from $C$ forming dot product $\lambda_j$ with $x$ is exactly~$\tau_j$.
\end {theorem}
Most of the proof of Theorem \ref {weakly2msharp} repeats the proof of Theorem \ref {stiff}. However, item (i) in Theorem \ref {weakly2msharp} requires a longer proof and the potential function in Theorem \ref {weakly2msharp} has a finite value at $t=1$ which must also be an interpolation point. Therefore, we do not combine the two proofs. This also helps with their readability.

\begin {proof}[Proof of Theorem \ref {weakly2msharp}]
Let $q\in \mathbb P_{2m+\nu}$ be the polynomial interpolating $f$ and $f'$ at points $\lambda_1,\ldots,\lambda_m$ which interpolates $f$ at $\lambda_{m+1}=1$, and if $\nu=1$, interpolates $f$ at $\lambda_0=-1$. For every $t\in [-1,1]$, there is a point $\xi\in (-1,1)$ such that 
\begin {equation}\label {c_1}
f(t)-q(t)=\frac {f^{(2m+1+\nu)}(\xi)}{(2m+1+\nu)!}(t-1)(t-\lambda_1)^2\cdots (t-\lambda_m)^2(t+1)^\nu.
\end {equation}
Since, by assumption, $f^{(2m+1+\nu)}\geq 0$ on $(-1,1)$, we have $f(t)-q(t)\leq 0$ or $q(t)\geq f(t)$, $t\in [-1,1]$. Since $(C,W)$ is a $(2m+\nu)$-design and quadrature \eqref {quadL} is exact on $q(t)$, for every point $x\in \mathbb{S}^{n-1}$, using \eqref {designequiv} we have
\begin {equation}\label {pq}
\begin {split}
U_f(x,C;W)&=\sum\limits_{i=1}^{N}w_if(x\cdot x_i)\leq \sum\limits_{i=1}^{N}w_iq(x\cdot x_i)=\alpha_0(q)\\
&=\sum\limits_{j=1-\nu}^{m+1}\tau_jq(\lambda_j)=\sum\limits_{j=1-\nu}^{m+1}\tau_jf(\lambda_j).
\end {split}
\end {equation}

Let now $(C,W)$ be a weakly sharp $(2m+\nu)$-design. Let also $z\in A_{1,\nu}$ be any point that forms $m+\nu$ distinct values of the dot product with other points from $C$. By Lemma \ref {allofthem}, the values of these dot products plus the dot product $1$ are the zeros $\lambda_{1-\nu}<\ldots<\lambda_{m+1}$ of the polynomial $(1-t)(1+t)^\nu P_m^{(1,\nu)}(t)$ and for each $j=1-\nu,\ldots,m+1$, the sum of the weights of points in $C$ that form dot product $\lambda_j$ with $z$ is $\tau_j$. Then
$$
U_f(z,C;W)=\sum\limits_{i=1}^{N}w_if(z\cdot x_i)=\sum\limits_{j=1-\nu}^{m+1}\tau_jf(\lambda_j).
$$

Assume now that $f^{(2m+1+\nu)}>0$ on $(-1,1)$. Then in view of \eqref {c_1} we have $f(t)<q(t)$ for $t\in [-1,1]\setminus \{\lambda_{1-\nu},\ldots,\lambda_{m+1}\}$. Let $x\in \mathbb{S}^{n-1}$ be any point such that equality holds in \eqref {potential-upper}. Then we have equality in the first line of \eqref {pq}. Since $f\leq q$ on $[-1,1]$, we cannot have $f(x\cdot  x_i)<q(x\cdot x_i)$ for any $i$ without ruining this equality. Therefore, each dot product $ x\cdot x_i$ is a point of equality of $f$ and $q$; that is, one of the numbers $\lambda_{1-\nu},\ldots,\lambda_{m+1}$.

If $\nu=0$, since $(C,W)$ is a $2m$-design, by Lemma~\ref {allofthem1} with $\alpha=\beta=0$, dot product $x\cdot x_i$ assumes at least $m+1$ distinct values. Then it assumes each of the values $\lambda_{1},\ldots,\lambda_{m+1}$. In particular, $x\in C=A_{1,0}$.

If $\nu=1$, since $(C,W)$ is a $2m$-design, by Lemma~\ref {allofthem1} with $\alpha=\beta=0$, the point $x$ must form at least $m+1$ dot products with points from $C$. Then $x$ must form at least one of the dot products $\lambda_0=-1$ or $\lambda_{m+1}=1$ with points from $C$; that is, $x\in C\cup (-C)$. We now apply again Lemma \ref {allofthem1} with $\alpha=1$ and $\beta=0$ if $x\in C=A_{1,0}$ or with $\alpha=0$ and $\beta=1$ if $x\in -C=A_{0,1}$. Since $(C,W)$ is a $(2m+\alpha+\beta)$-design, $x$ must form at least $m+2$ distinct dot products with points from $C$ thus forming each dot product $\lambda_0,\ldots,\lambda_{m+1}$. In particular, $x\in C\cap (-C)=A_{1,1}$.

Thus, in both cases, $x\in A_{1,\nu}$ and dot product $x\cdot x_i$ assumes exactly $m+1+\nu$ values, which are $\lambda_{1-\nu},\ldots,\lambda_{m+1}$; that is, $(C,W)$ is a weakly sharp $(2m+\nu)$-design. By Lemma \ref {allofthem}, the sum of the weights of points in $C$ forming dot product $\lambda_j$ with $x$ is $\tau_j$, $j=1-\nu,\ldots,m+1$. 
\end {proof}

\section {Examples}

In this section, we will first consider positive linear combinations of weighted designs of equal strength to obtain a weighted design of the same strength, see Proposition \ref {convexcomb}. This construction does not depend on the mutual placement of the designs on the sphere. \textcolor{black}{Utilizing a careful mutual placement of the component designs and their weights yields weighted designs of a higher strength. The corresponding bounds will be explored in a future work.}

Given two weighted spherical codes $(C_1,W_1)$ and $(C_2,W_2)$, and positive numbers $\alpha,\beta$ such that $\alpha+\beta=1$, denote by $\alpha (C_1,W_1)+\beta(C_2,W_2)$ the code $(C_1\cup C_2,W)$, where each point $x\in C_1\setminus C_2$, has weight $\alpha$ times the weight of $x$ in $(C_1,W_1)$, each point $x\in C_2\setminus C_1$ has weight $\beta$ times the weight of $x$ in $(C_2,W_2)$, and each point in $C_1\cap C_2$ has weight $\alpha w+\beta v$, where $w$ is the weight of $x$ in $(C_1,W_1)$ and $v$ is the weight of $x$ in $(C_2,W_2)$. We will call this code the weighted union of $(C_1,W_1)$ and $(C_2,W_2)$. It corresponds to the discrete probability measure $\kappa(C_1\cup C_2,W)=\alpha\cdot\kappa(C_1,W_1)+\beta\cdot\kappa(C_2,W_2)$.

\begin {proposition}\label {convexcomb}
Let $(C_1,W_1)$ and $(C_2,W_2)$ be spherical $m$-designs on $\mathbb{S}^{n-1}$, $n\geq 2$, $m\geq 1$, and $\alpha,\beta>0$ be real numbers such that $\alpha+\beta=1$. Then their weighted union $\alpha (C_1,W_1)+\beta(C_2,W_2)$ is also an $m$-design.
\end {proposition}
\begin {proof}
Let $w(z)$ be the weight of a point $z$ in $\alpha (C_1,W_1)+\beta(C_2,W_2)$, $w_1(z)$ be the weight of a point $z$ in $(C_1,W_1)$ and $w_2(z)$ be the weight of a point $z$ in $(C_2,W_2)$.
For every polynomial $p\in \mathbb P_m$ and every $x\in \mathbb{S}^{n-1}$,
\begin {equation*}
\begin {split}
\sum\limits_{z\in C_1\cup C_2}w(z)p(x\cdot z)&=\sum\limits_{z\in C_1\setminus C_2}\alpha w_1(z)p(x\cdot z)+\sum\limits_{z\in C_2\setminus C_1}\beta w_2(z)p(x\cdot z)\\
&\quad +\sum\limits_{z\in C_1\cap C_2}(\alpha w_1(z)+\beta w_2(z))p(x\cdot z)\\
&=\alpha \sum\limits_{z\in C_1}w_1(z)p(x\cdot z)+\beta \sum\limits_{z\in C_2}w_2(z)p(x\cdot z)\\
&=\alpha \alpha_0 (p)+\beta \alpha_0(p)=\alpha_0(p);
\end {split}
\end {equation*}
that is, the weighted union is also an $m$-design.
\end {proof}

An example of a weighted design attaining both universal bounds is the square-based pyramid we mentioned earlier. It is the weighted code $(C,W)\in J_{3,5}$, where $C\subset \mathbb{S}^2$ consists of the point $y=(0,0,1)$ with a weight $1/4$ and four points $(\pm 2/3,\pm 2/3,-1/3)$ each with a weight $3/16$. It is a weighted $2$-design with the point $y\in C$ forming only one dot product (which is $-1/3$) with other points of the code. Thus, it is a weakly sharp $2$-design.
The adjacent polynomial $P_1^{(1,0)}$ has a zero at $\gamma_1=-1/3$ and the corresponding weights of quadrature \eqref {quad1} are
$$
\sigma_1=\frac {1}{2}\int_{-1}^{1}R_1(t)\ \! dt=\frac {1}{2}\int_{-1}^{1}\frac {(1-t)}{4/3}\ \! dt=\frac {3}{4}
$$
and
$$
\sigma_2=\frac {1}{2}\int_{-1}^{1}R_2(t)\ \! dt=\frac {1}{2}\int_{-1}^{1}\(\frac {t+1/3}{4/3}\)^2\ \! dt=\frac {1}{4},
$$
where $R_j$, $j=1,2$, is the polynomial defined in \eqref {Rj} and \eqref {R0}. Then Theorem \ref {weakly2msharp} implies that for every admissible potential function $f:[-1,1]\to\RR$ such that $f'''\geq 0$ on $(-1,1)$, we have
$$
U_f(x,C;W)\leq \frac {1}{4}f(1)+\frac {3}{4}f\(-\frac {1}{3}\),\ \ \ x\in \mathbb{S}^2.
$$
This bound is sharp and is attained at $y=(0,0,1)$. If $f'''>0$ on $(-1,1)$, then this bound is attained only at $y$ because $y$ is the only point in the code forming one dot product with other points of the code. Similarly, by Theorem \ref {stiff}, for every $$
U_f(x,C;W)\geq \frac {1}{4}f(-1)+\frac {3}{4}f\(\frac {1}{3}\),\ \ \ x\in \mathbb{S}^2.
$$
This bound is sharp and is attained at $-y=(0,0,-1)$. If $f'''>0$ on $(-1,1)$, then this bound is attained only at $-y$ because $-y$ is the only point in the code $-C$ forming two dot products with points of the code one of which is $-1$.

A second example of a weighted design satisfying the assumptions of Theorem~\ref {stiff} and attaining the universal lower bound is the set $C\subset \mathbb{S}^{n-1}$, $n\geq 4$, of vertices of a hypercube obtained as a weighted union of two demihypercubes. Every point of $C$ is a vector of the form $\(\pm \frac {1}{\sqrt{n}},\ldots,\pm \frac {1}{\sqrt{n}}\)\in \RR^n$. Let $C_0$ be the set of all such vectors with an even number of minus signs and $C_1$ be the set of all such vectors with an odd number of minus signs. Each set $C_0$ and $C_1$ is an equi-weighted $3$-design, see \cite {Bor-FL}, of cardinality $2^{n-1}$. We assign a weight $w_0>0$ to each point of $C_0$ and a weight $w_1>0$ to each point of $C_1$ so that $w_0+w_1=2^{1-n}$. By Proposition \ref {convexcomb}, this yields a weighted $3$-design $(C,W)$. The set of points on $\mathbb{S}^{n-1}$ forming two distinct dot products with points of $C$ is the cross-polytope $\Omega_{2n}^\ast=\{\pm e_1,\ldots,\pm e_n\}$, where $e_1,\ldots,e_n$ are the standard basis vectors in $\RR^n$. Then $(C,W)$ is $2$-stiff. The Gegenbauer polynomial $P_2^{(n)}(t)=\frac {nt^2-1}{n-1}$ has zeros $\alpha_1=-\frac {1}{\sqrt{n}}$ and $\alpha_2=\frac {1}{\sqrt{n}}$. Each vector $\pm e_j$ from $\Omega_{2n}^\ast$ forms dot product $\frac {1}{\sqrt{n}}$ with exactly half of points of $C_0$ and with exactly half of points of $C_1$. Thus, the total weight of points from $C$ forming dot product $\frac {1}{\sqrt{n}}$ with $\pm e_j$ is $\frac {1}{2}$. Then the total weight of points from $C$ forming dot product $-\frac {1}{\sqrt{n}}$ with $\pm e_j$ is also $\frac {1}{2}$. By Theorem~\ref {stiff}, for any admissible potential function $f:[-1,1]\to (-\infty,\infty]$ with $f^{(4)}\geq 0$ on $(-1,1)$, the universal lower bound 
$$
U_f(x,C;W)\geq \frac {1}{2}\(f\(-\frac {1}{\sqrt{n}}\)+f\(\frac {1}{\sqrt{n}}\)\),\ \ \ x\in \mathbb{S}^{n-1},
$$ 
holds. It is sharp and is attained at every point $x\in \Omega_{2n}^\ast$. If, in addition, $f^{(4)}>0$ on $(-1,1)$, then it is attained only at points $x\in \Omega_{2n}^\ast$.

A third example of weighted designs that attain the universal lower bound can be constructed from any $m$-stiff configuration $\W C_N$ on $\mathbb{S}^{n-1}$ which does not have an antipodal pair. Such codes can be found on $\mathbb{S}^{20}$ for $N=112$ or $N=162$ and $m=2$ and on $\mathbb{S}^{21}$ for $N=100$ and $m=2$ (the Higman-Simms configuration) or for $N=891$ and $m=3$, see \cite {BDHSS-Sharp}. They also exist on $S^{2k}$, $k\geq 2$, with $N=2^{2k}$ and $m=2$ mentioned above as a demihypercubes in $\RR^n$ for $n\geq 5$ odd\footnote {For $n\geq 4$ even, the demihypercube on the unit sphere in $\RR^n$ is an antipodal configuration.}. 

We construct the weighted code as the set $C_{2N}:=\W C_{N}\cup (-\W C_{N})$ assigning weight $\alpha>0$ to each point of $C_{N}$ and weight $\beta>0$ to each point of $-C_{N}$ such that $\alpha+\beta=\frac {1}{N}$. By Proposition \ref {convexcomb}, the disjoint union $C_{2N}$ is a weighted $(2m-1)$-design. It is $m$-stiff, since by Lemma \ref {allofthem}, for every point $z\in \mathbb{S}^{n-1}$ forming $m$ distinct dot products with points of $\W C_N$, those dot products are zeros $-1<\alpha_1<\ldots<\alpha_m<1$ of the Gegenbauer polynomial $P_m^{(n)}$ and, hence, are symmetric about~$0$. Then $z$ forms the same set of dot products with points of $-C_N$. Thus, $z$ forms $m$ distinct dot products with points of $C_{2N}$. By Theorem~\ref {stiff}, for any admissible potential function $f:[-1,1]\to (-\infty,\infty]$ with $f^{(2m)}\geq 0$ on $(-1,1)$, the universal lower bound 
\begin {equation}\label {lower}
U_f(x,C_{2N};W)\geq \sum\limits_{j=1}^{m}\rho_jf(\alpha_j),\ \ \ x\in \mathbb{S}^{n-1},
\end {equation}
holds. It is sharp and is attained at every point $x\in \mathbb{S}^{n-1}$ forming $m$ distinct dot products with points of $C_{2N}$. If, in addition, $f^{(2m)}>0$ on $(-1,1)$, then bound \eqref {lower} is attained only at such points.



\section {Constrained frame energy minimization}\label {7c}

This section is devoted to energy minimization problem for a class of charges that contains certain discrete distributions and certain charges with a discrete and continuous (singular) component; the total charge is fixed. We start with a general result.
A Borel probability measure $\mu$ on $\mathbb{S}^{n-1}$ is called a $(k,k)$-design, $k\geq 1$, if 
$$
\int\limits_{\mathbb{S}^{n-1}}\int\limits_{\mathbb{S}^{n-1}}P_{2\ell}^{(n)}\(x\cdot y\)\ \! d\mu(x) d\mu(y)=0,\ \ \ \ell=1,\ldots,k.
$$
By default, any Borel probability measure on $\mathbb{S}^{n-1}$ is a $(0,0)$-design.
Given a weighted code $(C,W)=(x_1,\ldots,x_N;w_1,\ldots,w_N)$, recall that
$
\kappa(C,W)=\sum\limits_{i=1}^{N}w_i\delta_{x_i},
$
where $\delta_x$ is the Dirac delta-measure with mass $1$ at $x$. The weighted code $(C,W)$ is called a weighted $(k,k)$-design, $k\geq 1$, if the measure $\kappa(C,W)$ is a $(k,k)$-design; i.e., if
\begin {equation}\label {equiv1}
\sum\limits_{i=1}^{N}\sum\limits_{j=1}^{N}w_iw_jP_{2\ell}^{(n)}\(x_i\cdot x_j\)=0,\ \ \ \ell=1,\ldots,k.
\end {equation}
For a given set  $K\subset [-1,1]$, let 
$$
D_K:=\{(x,y)\in \mathbb{S}^{n-1}\times \mathbb{S}^{n-1} : x\cdot y\in K\}.
$$
Let $\mu$ be any Borel probability measure on $\mathbb{S}^{n-1}$, $n\geq 2$. Denote
$$
\theta_\mu:=(\mu\times\mu)(D_{\{-1,1\}}).
$$
If $\mu=\kappa(C,W)$, we have
\begin {equation}\label {theta_mu}
\theta_\mu=\sum\limits_{i=1}^Nw_i(w_i+w_i'),
\end {equation}
where $w_i'$ is the weight corresponding to $-x_i$ (if $-x_i\not\in C$, we let $w_i'=0$).
For any continuous function $f:[0,1]\to \RR$, we define the $f$-energy of $\mu$ as
\begin {equation}\label {energy2}
E_f(\mu):=\int\limits_{\mathbb{S}^{n-1}}\int\limits_{\mathbb{S}^{n-1}}f\((x\cdot y)^2\)\ \! d\mu(y)\ \! d\mu(x).
\end {equation}
For a weighted code $(C,W)$, its $f$-energy is the $f$-energy of the measure $\mu=\kappa (C,W)$:
$$
E_f(C,W):=E_f(\mu)=\sum\limits_{x,y\in C}\mu\{x\}\mu\{y\}f\((x\cdot y)^2\)=\sum\limits_{i,j=1}^{N}w_iw_jf\((x_i\cdot x_j)^2\).
$$
When $f(t)=\left|t\right|^{p/2}$, $p>0$, in \eqref {energy2}, we obtain the $p$-frame energy of a general Borel probability measure on the sphere. In particular, the $p$-{\it frame energy} of a weighted code $(C,W)$ equals
$$
E^p(C,W):=\sum\limits_{i,j=1}^{N}w_iw_j\left|x_i\cdot x_j\right|^p.
$$


The following general result holds. 
\begin {theorem}\label {genframe}
Let $n\geq 2$, $m\geq 1$, and $A$ be a set of numbers $-1<\alpha_1<\ldots<\alpha_m<1$ symmetric about the origin. Suppose there exists a weighted $(m-1,m-1)$-design $(C,W)$ on $\mathbb{S}^{n-1}$ such that any dot product formed by a pair of points from $C$ lies in the set $A\cup \{-1,1\}$.

Let $f:[0,1]\to \RR$ be a continuous function such that $f^{(k)}\geq 0$ in $(0,1)$, $k=1,\ldots,m$. 
If $\mu$ is any Borel probability measure with $\theta_\mu\geq \theta_{\mu^\ast}$, where $\mu^\ast=\kappa(C,W)$, then
\begin {equation}\label {en_ineq}
E_f(\mu)\geq E_f(\mu^\ast)=\sum\limits_{i,j=1}^{N}w_iw_jf\((x_i\cdot x_j)^2\).
\end {equation}

Assume, in addition, that $f^{(k)}>0$ on $(0,1)$ for $k=1,\ldots,m$. Then equality in \eqref {en_ineq} holds if and only if $\mu$ is a $(m-1,m-1)$-design, $\(\mu\times \mu\)\(D_{(-1,1)\setminus A}\)=0$, and $\theta_\mu=\theta_{\mu^\ast}$.

\end {theorem}
We recall that for potential functions with positive first $M$ derivatives and negative $(M+1)$-th derivative, Bilyk et. al. \cite {BilGlaPar2022} established the energy minimizing property of equi-weighted discrete measures supported at points of tight spherical $M$-designs among all Borel probability measures on~$\mathbb{S}^{n-1}$. 
The optimal configuration in Theorem \ref {genframe} need not be a tight design due to a weaker requirement for the strength of the optimal design.
The clustering phenomenon of measures minimizing the energy integral on the sphere was also studied in \cite {BilGlaPar2021}.

\begin {remark}
{\rm When $f(t)=t^{p/2}$ (we call $f$ the $p$-frame potential in this case), the assumption that $f^{(k)}>0$ in $(0,1)$ for $k=1,\ldots,m$ holds when $p> 2m-2$.}
\end {remark}

\begin {proof}[Proof of Theorem \ref {genframe}]
Denote by $L\geq 0$ and $\nu\in \{0,1\}$ integers such that $m=2L+\nu$ and let $\beta_i:=2\alpha_i^2-1$, $i=1,\ldots,L+\nu$. Denote by $p(t)$ the unique polynomial in $\mathbb P_{m-1}$ interpolating values of the function $g(t):=f\(\frac {t+1}{2}\)$, $t\in [-1,1]$, at points $\beta_i$, $i=1,\ldots,L+\nu$, and values of the derivative $g'$ at $\beta_i$, $i=1,\ldots,L$. For every $t\in [-1,1]$, there exists a point $\xi=\xi(t)\in (-1,1)$ such that 
\begin {equation}\label {Hermite}
g(t)-p(t)=\frac {g^{(m)}(\xi)}{m!}(t+1)^\nu(t-\beta_1)^2\cdots (t-\beta_L)^2.
\end {equation}
Since $g^{(m)}\geq 0$, $t\in (-1,1)$, we have $g(t)\geq p(t)$, $t\in [-1,1]$. Let
$$
p(2t^2-1)=d_0+\sum\limits_{k=1}^{m-1}d_{k}P_{2k}^{(n)}(t)
$$
be the Gegenbauer expansion of the even polynomial $p(2t^2-1)$.
Taking into account the fact that $f(t)=g(2t-1)\geq p(2t-1)$, $t\in [0,1]$, we have
\begin {equation}\label {interp3}
\begin {split}
E_f(\mu)&=\iint\limits_{D_{(-1,1)}}f\((x\cdot y)^2\)\ \! d(\mu\times\mu)(x,y)+\iint\limits_{D_{\{-1,1\}}}f\((x\cdot y)^2\)\ \! d(\mu\times\mu)(x,y)\\
&\geq G(\mu):=\iint\limits_{D_{(-1,1)}}p\(2(x\cdot y)^2-1\)\ \! d(\mu\times\mu)(x,y)+f(1)(\mu\times\mu)\(D_{\{-1,1\}}\)\\
&=\int\limits_{\mathbb{S}^{n-1}\times \mathbb{S}^{n-1}}p\(2(x\cdot y)^2-1\)\ \! d(\mu\times\mu)(x,y)+(f(1)-p(1))(\mu\times\mu)\(D_{\{-1,1\}}\)\\
&=d_0+\sum\limits_{k=1}^{m-1}d_{k}\int\limits_{\mathbb{S}^{n-1}}\int\limits_{\mathbb{S}^{n-1}}P_{2k}^{(n)}\(x\cdot y\)\ \! d\mu(y)\!\ d\mu(x)+(f(1)-p(1))\theta_\mu.
\end {split}
\end {equation}
The classical positive definiteness property of Gegenbauer polynomials implies that
\begin {equation}\label {GP}
\int\limits_{\mathbb{S}^{n-1}}\int\limits_{\mathbb{S}^{n-1}}P_{2k}^{(n)}\(x\cdot y\)\ \! d\mu(y)\!\ d\mu(x)\geq 0,\ \ \ k=1,\ldots,m-1.
\end {equation}
We have $f(1)\geq p(1)$ and, by assumption, $\theta_\mu\geq \theta_{\mu^\ast}$. By Theorem \ref {posdefeven} in Appendix, we have $d_{k}\geq 0$, $k=1,\ldots,m-1$, when $m\geq 2$, and, consequently,
\begin {equation}\label {lbound}
E_f(\mu)\geq G(\mu)\geq d_0+(f(1)-p(1))\theta_\mu\geq d_0+(f(1)-p(1))\theta_{\mu^\ast}.
\end {equation}
For $m=1$, estimate \eqref {lbound} holds without using Theorem \ref {posdefeven}. At the same time, for the weighted $(m-1,m-1)$-design $(C,W)$, we have
\begin {equation*}
\begin {split}
E_f(\mu^\ast)&=\sum\limits_{i=1}^{N}\sum\limits_{j=1}^{N}w_iw_j f\((x_i\cdot x_j)^2\)\\
&=\sum\limits_{x\in C}\sum\limits_{y\in C\atop y\neq x,-x}\mu^\ast\{x\}\mu^\ast\{y\}f\((x\cdot y)^2\)+f(1)\sum\limits_{x\in C}\mu^\ast\{x\}\mu^\ast\{x,-x\}.
\end {split}
\end {equation*}
For any $y\neq x,-x$, dot product $x\cdot y$ equals $\alpha_j$ for some $1\leq j\leq m$. Since $p$ interpolates $g$ at each $\beta_i=2\alpha_i^2-1$, $i=1,\ldots,L+\nu$, we have $f\((x\cdot y)^2\)=g\(2\alpha_j^2-1\)=g(\beta_j)=p(\beta_j)=p(2(x\cdot y)^2-1)$. Since $(C,W)$ is a weighted $(m-1,m-1)$-design, we have
\begin {equation*}
\begin {split}
E_f(\mu^\ast)&=\sum\limits_{x\in C}\sum\limits_{y\in C\atop y\neq x,-x}\mu^\ast\{x\}\mu^\ast\{y\}p\(2(x\cdot y)^2-1\)+f(1)\sum\limits_{x\in C}\mu^\ast\{x\}\mu^\ast\{x,-x\}\\
&=\sum\limits_{x\in C}\sum\limits_{y\in C}\mu^\ast\{x\}\mu^\ast\{y\}p\(2(x\cdot y)^2-1\)+(f(1)-p(1))\sum\limits_{x\in C}\mu^\ast\{x\}\mu^\ast\{x,-x\}\\
&=d_0+\sum\limits_{k=1}^{m-1}d_{k}\sum\limits_{i,j=1}^{N}w_iw_jP_{2k}^{(n)}(x_i\cdot x_j)+(f(1)-p(1))\theta_{\mu^\ast}\\
&=d_0+(f(1)-p(1))\theta_{\mu^\ast}.
\end {split}
\end {equation*}
Combined with \eqref {lbound}, this yields \eqref {en_ineq}.

Assume now that $f^{(k)}>0$ on $(0,1)$, $k=1,\ldots,m$. Then $g^{(k)}(t)>0$, $t\in (-1,1)$, $k=1,\ldots,m$. By \eqref {Hermite}, we have $g(t)>p(t)$ for all $t\in [-1,1]\setminus \{\beta_1,\ldots,\beta_{L+\nu}\}$ and $g(t)=p(t)$, $t\in \{\beta_1,\ldots,\beta_{L+\nu}\}$. Consequently, $f\((x\cdot y)^2\)=g(2(x\cdot y)^2-1)>p\(2(x\cdot y)^2-1\)$, $(x,y)\in D_{(-1,1)\setminus A}=D_{(-1,1)}\setminus D_A$ and $f\((x\cdot y)^2\)=g(2(x\cdot y)^2-1)=p\(2(x\cdot y)^2-1\)$, $(x,y)\in D_{A}$. Then equality in \eqref {interp3} holds if and only if $\(\mu\times \mu\)\(D_{(-1,1)\setminus A}\)=0$. 

If $m\geq 2$, by Theorem \ref {posdefeven}, we have $d_{k}>0$, $k=1,\ldots,m-1$. The second inequality in \eqref {lbound} becomes an equality if and only if all integrals in \eqref {GP} equal zero; i.e., if and only if $\mu$ is a $(m-1,m-1)$-design with $\theta_{\mu}\geq \theta_{\mu^\ast}$. If $m=1$, then the second inequality in \eqref {GP} is an equality for any Borel probability measure $\mu$ with $\theta_\mu\geq \theta_{\mu^\ast}$. Since $\mu$ is a $(0,0)$-design by default, the second equality holds if and only if $\mu$  is a $(0,0)$-design with $\theta_\mu\geq \theta_{\mu^\ast}$. 

Since by \eqref {Hermite} $f(1)=g(1)>p(1)$, the third inequality in \eqref {lbound} becomes an equality if and only if $\theta_\mu=\theta_{\mu^\ast}$.
\end {proof}

 We will call a weighted code $(C,W)$ {\it antipodal} if for any point $x\in C$, the point $-x$ is in $C$ and the weights of $x$ and $-x$ are equal. 

\begin {remark}
{\rm
Recall that the real projective space $\mathbb R\mathbb P^{n-1}$, $n\geq 2$, can be considered as a metric space where points are pairs of antipodal vectors on $\mathbb{S}^{n-1}$ with the distance $d({\bf x},{\bf y})$ between pairs ${\bf x}$ and ${\bf y}$ represented by vectors $x$ and $y$ being defined by $\cos d({\bf x},{\bf y})=2(x\cdot y)^2-1$. Observe that $f\((\pm x\cdot (\pm y))^2\)=f\(\frac {1+\cos d({\bf x},{\bf y})}{2}\)=g\(\cos d({\bf x},{\bf y})\)$, where the function $g(t)=f\(\frac {1+t}{2}\)$ is used in the proof of Theorem \ref {genframe}. Then the $f$-energy of any weighted code 
$$
(C,W)=(x_1,\ldots,x_N,-x_1,\ldots,-x_N;w_1,\ldots,w_N,v_1,\ldots,v_N)
$$
consisting of $N$ distinct antipodal pairs with points in each pair having arbitrary non-negative weights with the sum of all $2N$ weights being $1$ can be written as the energy of the corresponding weighted code $\{{\bf x}_1,\ldots,{\bf x}_N;w_{{\bf x}_1},\ldots,w_{{\bf x}_N}\}$ in $\mathbb R\mathbb P^{n-1}$:
$$
E_f(C,W)=\sum\limits_{i,j=1}^{N}(w_i+v_i)(w_j+v_j)f\(\left|x_i\cdot x_j\right|^2\)=\sum\limits_{i.j=1}^N w_{{\bf x}_i}w_{{\bf x}_j}g\(\cos d({\bf x}_i,{\bf x}_j)\),
$$
where each ${\bf x}_i$ is represented by vectors $x_i$ and $-x_i$ and $w_{{\bf x}_i}=w_i+v_i$, $i=1,\ldots,N$. This also means that the energy $E_f(C,W)$ of a weighted code $(C,W)$ with an antipodal code $C$ is independent of the distribution of the weights within each antipodal pair of $C$.
}
\end {remark}
Theorem \ref {genframe} has the following consequence for weighted codes. In view of the above remark, Theorem \ref {antipodal} has an analogue in the projective space. We recall that the universal energy minimization property of sharp codes among equi-weighted point configurations was established by Cohn and Kumar in their seminal paper \cite {CohKum2007} for any two-point homogeneous metric space which includes the Euclidean sphere and the real projective space. 
\begin {theorem}\label {antipodal}
Let $n\geq 2$, $m,N\geq 1$, and $A$ be a set of numbers $-1<\alpha_1<\ldots<\alpha_m<1$ symmetric about the origin. Suppose there exists an equi-weighted $N$-point $(m-1,m-1)$-design $C^\ast=\{x_1^\ast,\ldots,x_N^\ast\}\subset \mathbb{S}^{n-1}$ with no antipodal pairs and $x_i^\ast\cdot x_j^\ast\in A$, $1\leq i\neq j\leq N$. 
Let $f:[0,1]\to \RR$ be a continuous function such that $f^{(k)}\geq 0$ in $(0,1)$ for $k=1,\ldots,m$. 
Suppose $(C,W)=(x_1,\ldots,x_{2N};w_1,\ldots,w_{2N})$, where $C$ consists of $N$ distinct antipodal pairs on $\mathbb{S}^{n-1}$ and $W=(w_1,\ldots,w_{2N})$ is any vector of non-negative numbers that add up to $1$. Then
\begin {equation}\label {regular}
E_f(C,W)=\sum\limits_{i,j=1}^{2N}w_iw_jf\((x_i\cdot x_j)^2\)\geq E_f(C^\ast,W^\ast)=\frac {1}{N^2}\sum\limits_{i,j=1}^{N}f\((x_i^\ast\cdot x_j^\ast)^2\),
\end {equation}
where $W^\ast=\(\frac {1}{N},\ldots,\frac {1}{N}\)\in \RR^N$.

Assume, in addition, that $f^{(k)}>0$ on $(0,1)$ for $k=1,\ldots,m$. 
Then 
equality in \eqref {regular} holds if and only if $(C,W)$ is a weighted $(m-1,m-1)$-design on $\mathbb{S}^{n-1}$, where each antipodal pair in $C$ has a total weight of $\frac {1}{N}$, and dot products between points from distinct antipodal pairs lie in $A$.
\end {theorem}

\begin {proof}
Select one point from each antipodal pair of the code $C$ and denote the selected points by $y_1,\ldots,y_N$. Then $C=\{\pm y_1,\ldots,\pm y_N\}$. Letting $\mu=\sum\limits_{i=1}^{2N}w_i\delta_{x_i}$ and $\epsilon_i=\mu\{y_i,-y_i\}-\frac {1}{N}$, $i=1,\ldots,N$, we have $\sum_{i=1}^{N}\epsilon_i=0$ and
\begin {equation}\label {thetamu}
\begin {split}
\theta_\mu &=\sum\limits_{i=1}^{N}\(\mu\{y_i\}\mu\{y_i,-y_i\}+\mu\{-y_i\}\mu\{y_i,-y_i\}\)=\sum\limits_{i=1}^{N}\(\mu\{y_i,-y_i\}\)^2\\
&=\sum\limits_{i=1}^{N}\(\epsilon_i+\frac {1}{N}\)^2=\sum\limits_{i=1}^{N}\epsilon_i^2+\frac {2}{N}\sum\limits_{i=1}^{N}\epsilon_i+\frac {1}{N}= \sum\limits_{i=1}^{N}\epsilon_i^2+\frac {1}{N}\geq \frac {1}{N}=\theta_{\mu^\ast},\nonumber
\end {split}
\end {equation}
where $\mu^\ast=\frac {1}{N}\sum_{i=1}^{N}\delta_{x_i^\ast}$. Applying Theorem \ref 
{genframe}, we have
\begin {equation}\label {mu_star}
E_f(C,W)=E_f(\mu)=\sum\limits_{i,j=1}^{2N}w_iw_jf\((x_i\cdot x_j)^2\)\geq E_f(\mu^\ast)=\frac {1}{N^2}\sum\limits_{i,j=1}^{N}f\((x_i^\ast\cdot x_j^\ast)^2\)=E_f(C^\ast,W^\ast).
\end {equation}
Assume now that $f^{(k)}>0$ on $(0,1)$ for $k=1,\ldots,m$.
By Theorem \ref {genframe}, equality holds in \eqref {mu_star} for a given $(C,W)$ if and only if $(C,W)$ is a weighted $(m-1,m-1)$-design,  $\(\mu\times \mu\)\(D_{(-1,1)\setminus A}\)=0$, and $\theta_\mu=\theta_{\mu^\ast}$. We have $\theta_\mu=\theta_{\mu^\ast}$ if and only if $\epsilon_i=0$, $i=1,\ldots,N$; that is, $\mu\{y_i,-y_i\}=\frac {1}{N}$, $i=1,\ldots,N$.

Let $(C,W)$ be such that equality holds in \eqref {mu_star}. Take any two points, $z_1$ and $z_2$ in $C$ from distinct antipodal pairs. Then $\mu\{z_i,-z_i\}=\frac {1}{N}$, $i=1,2$, and we can guarantee that $\mu\{z_i\}>0$, $i=1,2$, replacing $z_i$ with its antipode if necessary. If it were that $z_1\cdot z_2\notin A$, then the pair $(z_1,z_2)$ would be in $D_{(-1,1)\setminus A}$. Since $(\mu\times \mu)\{(z_1,z_2)\}=\mu\{z_1\}\mu\{z_2\}>0$, the condition $\(\mu\times \mu\)\(D_{(-1,1)\setminus A}\)=0$ would not hold. Therefore, $z_1\cdot z_2\in A$. Thus, $(C,W)$ is a weighted $(m-1,m-1)$-design on $\mathbb{S}^{n-1}$ with each antipodal pair having a total weight of $\frac {1}{N}$, and dot products between points from distinct antipodal pairs lying in $A$.
Conversely, if $(C,W)$ is a weighted $(m-1,m-1)$-design on $\mathbb{S}^{n-1}$ with each antipodal pair in $C$ having a total weight of $\frac {1}{N}$ and dot products between points from distinct antipodal pairs lying in $A$, then $\theta_\mu=\theta_{\mu^\ast}$.
Since for any $z_1,z_2\in C$ from distinct antipodal pairs, we have $z_1\cdot z_2\in A$, and for $z_1,z_2$ from the same antipodal pair, $z_1\cdot z_2=\pm 1$, the set $C\times C$ (which has $(\mu\times \mu)$-measure $1$) is contained in $D_{A\cup \{-1,1\}}$. Then $D_{A\cup \{-1,1\}}=\(\mathbb{S}^{n-1}\times \mathbb{S}^{n-1}\)\setminus D_{(-1,1)\setminus A}$ also has $(\mu\times \mu)$-measure $1$. Consequently, $\(\mu\times \mu\)\(D_{(-1,1)\setminus A}\)=0$. Then equality holds in \eqref {mu_star} and, hence, in \eqref {regular}.
\end {proof}

We now state the first consequence of Theorem \ref {antipodal}.
\begin {corollary}\label {nonantipodal}
Let $n\geq 2$, $m\geq 1$, $N\geq 2$, and $A$ be a set of numbers $-1<\alpha_1<\ldots<\alpha_m<1$ symmetric about the origin. Suppose there exists an equi-weighted $N$-point $(m-1,m-1)$-design $C^\ast=\{x_1^\ast,\ldots,x_N^\ast\}\subset \mathbb{S}^{n-1}$ with no antipodal pairs and $x_i^\ast\cdot x_j^\ast\in A$, $1\leq i\neq j\leq N$. 
Let $f:[0,1]\to \RR$ be a continuous function such that $f^{(k)}\geq 0$ in $(0,1)$ for $k=1,\ldots,m$. 
Suppose $(C,W)=(x_1,\ldots,x_N;w_1,\ldots,w_N)$ is any weighted code of $N$ distinct points on $\mathbb{S}^{n-1}$ and $W=(w_1,\ldots,w_N)$ is any vector of non-negative numbers that add up to $1$. Then
\begin {equation}\label {Nnonantipodal}
E_f(C,W)=\sum\limits_{i,j=1}^{N}w_iw_jf\((x_i\cdot x_j)^2\)\geq E_f(C^\ast,W^\ast)=\frac {1}{N^2}\sum\limits_{i,j=1}^{N}f\((x_i^\ast\cdot x_j^\ast)^2\),
\end {equation}
where $W^\ast=\(\frac {1}{N},\ldots,\frac {1}{N}\)\in \RR^N$.

Assume, in addition, that $f^{(k)}>0$ on $(0,1)$ for $k=1,\ldots,m$. Then equality holds in \eqref {Nnonantipodal} if and only if $(C,W)$ is an equi-weighted $(m-1,m-1)$-design with no antipodal pairs and dot products between distinct points lying in $A$.
\end {corollary}
\begin {proof}
Let $(C,W)$ be any weighted code of $N$ distinct points on $\mathbb{S}^{n-1}$. For each point in $C$ that does not have an antipode in $C$, we augment $(C,W)$ with its antipode with weight $0$. Then, if necessary, we add to $C$ antipodal pairs with total weight $0$ to get exactly $N$ antipodal pairs in $C$. Since we added to $C$ only points of weight $0$, the energy of the augmented weighted code $(C',W')$ satisfies $E_f(C',W')=E_f(C,W)$.  
Then by Theorem \ref {antipodal}, 
$$
E_f(C,W)=E_f(C',W')\geq E_f(C^\ast,W^\ast).
$$
Assume that $f^{(k)}>0$ on $(0,1)$ for $k=1,\ldots,m$ and that equality holds in \eqref {Nnonantipodal}.
 By Theorem~\ref {antipodal}, each antipodal pair in $C'$ has weight $\frac {1}{N}$. Therefore, we did not add to $C$ any antipodal pairs of total weight $0$. This also means that all $N$ points of $C$ do not have an antipode in $C$. Then the antipode of each point in $C$ has weight $0$ in $(C',W')$. Since each antipodal pair of $C'$ has a total weight $\frac {1}{N}$, each point of the original code $C$ must have a weight $\frac {1}{N}$; that is, $C$ is an equi-weighted $(m-1,m-1)$-design with no antipodal pairs. Theorem \ref {antipodal} also implies that dot products between distinct points in $C$ lie in $A$.

Conversely, every $N$-point equi-weighted $(m-1,m-1)$-design $(C,W)$ with no antipodal pairs and dot products between distinct points lying in $A$ can be augmented with antipodes of all its points with weight $0$ to obtain the antipodal code $C'$ and the vector $W'$ of $2N$ weights. Then the corresponding probability measure $\mu=\kappa(C',W')$ satisfies
$$
\theta_\mu=\sum\limits_{y\in C}\frac {1}{N}\(\frac {1}{N}+0\)+\sum\limits_{y\in -C}0\(0+\frac {1}{N}\)=\frac {1}{N}=\theta_{\mu^\ast}.
$$
The weighted code $(C',W')$ is a weighted $(m-1,m-1)$-design with each antipodal pair having a total weight $\frac {1}{N}$. Any two points $z_1,z_2\in C'$ from distinct antipodal pairs are modulo the sign two distinct points $y_1,y_2\in C$. Since $y_1\cdot y_2\in A$ and $A$ is symmetric, $z_1\cdot z_2\in A$. Then  
$E_f(C',W')=E_f(C^\ast,W^\ast)$ by Theorem \ref {antipodal}. Consequently, $E_f(C,W)=E_f(C',W')=E_f(C^\ast,W^\ast)$.
\end {proof}

We next state the second consequence of Theorem \ref {antipodal}.
\begin {corollary}\label {antipodal2N}
Let $n\geq 2$, $m,N\geq 1$, and $A$ be a set of numbers $-1<\alpha_1<\ldots<\alpha_m<1$ symmetric about the origin. Suppose there exists an antipodal equi-weighted $2N$-point $(2m-1)$-design $\W C=\{\W x_1,\ldots,\W x_N,-\W x_1,\ldots,-\W x_N\}\subset \mathbb{S}^{n-1}$ with dot products between points from distinct antipodal pairs lying in $A$.
Let $f:[0,1]\to \RR$ be a continuous function such that $f^{(k)}\geq 0$ in $(0,1)$ for $k=1,\ldots,m$. 
Suppose $(C,W)=(x_1,\ldots,x_{N},-x_1,\ldots,-x_N;w_1,\ldots,w_{N},w_1,\ldots,w_N)$ is any antipodal weighted code of $2N$ distinct points. Then
\begin {equation}\label {2N'}
E_f(C,W)=4\sum\limits_{i,j=1}^{N}w_iw_jf\((x_i\cdot x_j)^2\)\geq E_f(\W C,\W W)=\frac {1}{N^2}\sum\limits_{i,j=1}^{N}f\((\W x_i\cdot \W x_j)^2\),
\end {equation}
where $\W W=\(\frac {1}{2N},\ldots,\frac {1}{2N}\)\in \RR^{2N}$.

Assume, in addition, that $f^{(k)}>0$ on $(0,1)$ for $k=1,\ldots,m$. 
Then equality holds in \eqref {2N'} if and only if $(C,W)$ is an antipodal equi-weighted $2N$-point $(2m-1)$-design on $\mathbb{S}^{n-1}$ with dot products between points from distinct antipodal pairs lying in $A$.
\end {corollary}

The case $N=1$ of Corollary \ref {antipodal2N} is trivial, since the set of weighted antipodal $2$-point codes consists only of equi-weighted antipodal pairs, which all have the same energy (one should let $m=1$). The only other values of $N$ the case $m=1$ of Corollary~\ref {antipodal2N} applies to are $2\leq N\leq n$. The optimum is any equi-weighted code consisting of $N$ pairwise orthogonal antipodal pairs; i.e., a regular cross-polytope on $\mathbb{S}^{n-1}$ of dimension $N\leq n$.
 The set of vertices of an $N$-dimensional cross-polytope on $\mathbb{S}^{n-1}$ minimizes the $p$-frame energy among weighted antipodal $2N$-point codes for any $p>0$. The set of vertices of a regular simplex on $\mathbb{S}^{n-1}$, $n\geq 3$, symmetrized about the origin (it becomes a cube when $n=3$) minimizes the $p$-frame energy among weighted antipodal $N$-point codes for $N=2n+2$ and any $p\geq 2$. The set of vertices of the $24$-cell on $\mathbb{S}^3$ does so for $N=24$ and any $p\geq 4$.

\begin {remark}
{\rm
We call a code $C\subset \mathbb{S}^{n-1}$ {\it half-$m$-sharp} if for some $m\geq 1$, $C$ is an equi-weighted $(m-1,m-1)$-design on $\mathbb{S}^{n-1}$ with no antipodal pairs and distinct points of $C$ forming dot products that lie in a $m$-point set symmetric about $0$. If, in addition, $C$ is an equi-weighted $(m,m)$-design, then $C\cup (-C)$ is an $(m+1)$-sharp antipodal code; i.e., it is an antipodal $(m+1)$-distance set on $\mathbb{S}^{n-1}$, which is a $(2m+1)$-design. Such a code $C\cup (-C)$ will be a tight design, see \cite {DGS}.
For dimension $n$ and cardinality $N$ such that a half-sharp code $C$ of $N$ points exists on $\mathbb{S}^{n-1}$, by Corollary~\ref {nonantipodal}, the code $C$ is a solution to the $f$-energy minimization problem for arbitrary $N$-point weighted codes on $\mathbb{S}^{n-1}$ and by Corollary \ref {antipodal2N}, the code $C\cup (-C)$ is a solution to the $f$-energy minimization problem for antipodal $2N$-point weighted codes on $\mathbb{S}^{n-1}$.
}
\end {remark}

\begin {proof}[Proof of Corollary \ref {antipodal2N}]
By deleting one point from each antipodal pair of $(\W C,\W W)$ one can obtain an $N$-point equi-weighted $(m-1,m-1)$-design~$C^\ast$ with no antipodal pairs and dot products between distinct points lying in $A$. Let $W^\ast=\(\frac {1}{N},\ldots,\frac {1}{N}\)\in \RR^N$. Then Theorem \ref {antipodal} can be applied and for any weighted code $(C,W)=(x_1,\ldots,x_{N},-x_1,\ldots,-x_N;w_1,\ldots,w_{N},w_1,\ldots,w_N)$, we obtain that 
\begin {equation}\label {Ctilde}
E_f(C,W)\geq E_f(C^\ast,W^\ast)=E_f(\W C,\W W). 
\end {equation}
Assume that that $f^{(k)}>0$ on $(0,1)$ for $k=1,\ldots,m$. 
Let $(C,W)$ be any antipodal weighted code on $\mathbb{S}^{n-1}$ of $2N$ distinct points such that equality holds in \eqref {2N'}. Then equality holds in \eqref {regular}, and by Theorem \ref {antipodal}, $(C,W)$ is a weighted $(m-1,m-1)$-design with each antipodal pair having weight $\frac {1}{N}$ and dot products between points from distinct antipodal pairs lying in $A$. Since the weighted code $(C,W)$ is antipodal, it is a weighted $(2m-1)$-design. Since each antipodal pair must have weight $\frac {1}{N}$, each point of $C$ must have a weight of $\frac {1}{2N}$; i.e., $(C,W)$ is equi-weighted with dot products between points of $C$ from distinct antipodal pairs lying in $A$. 

Conversely, every antipodal equi-weighted $(2m-1)$-design $(C,W)$ on $\mathbb{S}^{n-1}$ with dot products between points from distinct antipodal pairs lying in $A$ is an $(m-1,m-1)$-design with each antipodal pair having a total weight $\frac {1}{N}$. By Theorem \ref {antipodal} it attains equality in \eqref {Ctilde} and, hence, in \eqref {2N'}.
\end {proof}

We conclude this section with an example for a regular simplex which follows from Corollary \ref {nonantipodal}.
\begin {proposition}
Let $n\geq 2$ and $(C^\ast,W^\ast)=(x_0^\ast,\ldots,x_n^\ast;\frac {1}{n+1},\ldots,\frac {1}{n+1})$ be an equi-weighted code on $\mathbb{S}^{n-1}$ with $\left|x_i^\ast\cdot x_j^\ast\right|=\frac {1}{n}$, $0\leq i\neq j\leq n$ (this includes the case of a regular $n$-simplex inscribed in $\mathbb{S}^{n-1}$). Suppose $f:[0,1]\to \RR$ is a continuous function such that $f',f''\geq 0$ on $(0,1)$. If $(C,W)=(x_0,\ldots,x_n;w_0,\ldots,w_n)$ is an arbitrary weighted $(n+1)$-point code on $\mathbb{S}^{n-1}$, where $W=(w_0,\ldots,w_n)$ is a vector of non-negative numbers that add up to $1$, then
\begin {equation}\label {simplexa}
\sum\limits_{i,j=0}^{n}w_iw_jf\((x_i\cdot x_j)^2\)\geq \frac {f(1)+nf\(\frac {1}{n^2}\)}{n+1}=\frac {1}{(n+1)^2}\sum\limits_{i,j=0}^{n}f\((x_i^\ast\cdot x_j^\ast)^2\).
\end {equation}
Assume, in addition, that $f',f''>0$ on $(0,1)$. Then equality in \eqref {simplexa} holds if and only if $w_0=\ldots=w_n=\frac {1}{n+1}$ and $\left|x_i\cdot x_j\right|=\frac {1}{n}$, $i\neq j$.
\end {proposition}
\begin {proof}
The equi-weighted code $(C^\ast,W^\ast)$ is a $(1,1)$-design. This is due to the fact that 
$$
\sum\limits_{i,j=0}^{n}P_2^{(n)}(x_i^\ast\cdot x_j^\ast)=(n+1)\(P_2^{(n)}(1)+nP_2^{(n)}\(\pm \frac {1}{n}\)\)=0.
$$
The regular $n$-simplex is an equi-weighted $(1,1)$-design with no antipodal pairs and dot products between distinct points lying in the set $A=\{-\frac {1}{n},\frac {1}{n}\}$. By Corollary \ref {nonantipodal}, we obtain \eqref {simplexa}.

Assume now that $f',f''>0$ on $(0,1)$. Assume that equality holds in \eqref {simplexa}. Then by Corollary \ref {nonantipodal}, $(C,W)$ is an equi-weighted $(n+1)$-point code with dot products between points with distinct indices being $\frac {1}{n}$ or $-\frac {1}{n}$. Conversely, if $(C,W)$ is an equi-weighted $(n+1)$-point code on $\mathbb{S}^{n-1}$ whose points with distinct indices form dot products $\frac {1}{n}$ or $-\frac {1}{n}$, then one can show directly the equality in \eqref {simplexa}. 
\end {proof}

\appendix

\section {Proof of a lemma for Theorem \ref {genframe}}

For completeness, we give the proof due to Levenshtein \cite {Lev} and Cohn and Kumar \cite {CohKum2007} of the positive definiteness of the Hermite interpolating polynomial for the potential function.
We start by recalling the following general result from~\cite {CohKum2007}.
\begin {theorem}\label {CKmain}
Let $\mu$ be any Borel measure on $\RR$ such that all polynomials are integrable with respect to $\mu$ with squares of non-zero polynomials having positive integrals over $\mu$. Suppose $\{p_k\}_{k=0}^{\infty}$ is the sequence of monic orthogonal polynomials with respect to $\mu$ such that ${\rm deg}\ \! p_k=k$, $k\geq 0$. Let $\gamma$ be any real number and let $r_1<r_2<\ldots<r_\ell$ be roots of $p_\ell(t)+\gamma p_{\ell-1}(t)$. 

Then for every $1\leq k<\ell$,
$$
\prod\limits_{i=1}^{k}(t-r_i)
$$
has positive coefficients in terms of $p_0(t),p_1(t),\ldots,p_k(t)$.
\end {theorem}

We will also need the following result from 
\cite {G}. 

\begin {theorem}\label {positiv}
If $\alpha\geq \beta>-1$, $\alpha+\beta+1\geq 0$, and $\{Q_k\}_{k=0}^{\infty}$ is the sequence of Jacobi polynomials orthogonal with weight $(1-t)^\alpha (1+t)^\beta$, then 
$$
\int\limits_{-1}^{1}Q_i(t)Q_j(t)Q_k(t)(1-t)^\alpha (1+t)^\beta \ \! dt\geq 0
$$
for all indices $i,j,k\geq 0$.
\end {theorem}

Given $n\geq 2$ and $\mu,\nu\in \{0,1\}$, denote by $\{P_k^{\mu,\nu}\}_{k=0}^{\infty}$ the sequence of monic polynomials orthogonal with the weight 
$$
v_{\mu,\nu}(t)=v_{\mu,\nu}^n(t):=(1+t)^{\nu-\frac {1}{2}}(1-t)^{\mu+\frac {n-3}{2}},\ \ \ t\in [-1,1],
$$
such that ${\rm deg}\  P^{\mu,\nu}_k=k$, $k\geq 0$, and let $P_k:=P_k^{0,0}$. 

\begin {lemma}\label {beta}
Let $n\geq 2$, $L\geq 2$, $\nu\in \{0,1\}$, and numbers $-1<\beta_L<\ldots<\beta_1<1$ be such that the polynomial  
$$
\Pi_L(t):=\prod\limits_{i=1}^{L}(t-\beta_i)
$$
is orthogonal to $\mathbb P_{L-2}$ with weight $v_{1,\nu}(t)$.
Then 
for $k=2,\ldots,L$,
$$
\prod\limits_{i=k}^{L}(t-\beta_i)=\sum\limits_{j=0}^{L-k+1}c_jP^{1,\nu}_j(t),\ \ \ \text{and}\ \ \ \prod\limits_{i=2}^{L}(t-\beta_i)\prod\limits_{i=k}^{L}(t-\beta_i)=\sum\limits_{j=0}^{2L-k}d_jP^{1,\nu}_j(t),
$$
where $\beta_i:=2\alpha_i^2-1$, $i=1,\ldots,L$, $c_j\geq 0$, $j=0,\ldots,L-k+1$, and $d_j\geq 0$, $j=0,\ldots,2L-k$.
\end {lemma}
\begin {proof}
Since $\Pi_L$ is orthogonal to $\mathbb P_{L-2}$ with weight $v_{1,\nu}(t)$, there exists $\gamma\in \RR$ such that $\Pi_L(t)=P_L^{1,\nu}(t)+\gamma P_{L-1}^{1,\nu}(t)$.
Theorem \ref {CKmain} now implies that for $k=2,\ldots,L$,
$$
\prod\limits_{i=k}^{L}(t-\beta_i)=\sum\limits_{j=0}^{L-k+1}c_{j,k}P^{1,\nu}_j(t),
$$
where $c_{j,k}\geq 0$, $j=0,\ldots,L-k+1$, $k=2,\ldots,L$. Furthermore,
$$
\prod\limits_{i=2}^{L}(t-\beta_i)\prod\limits_{i=k}^{L}(t-\beta_i)=\(\sum\limits_{j=0}^{L-1}c_{j,2}P^{1,\nu}_j(t)\)\(\sum\limits_{j=0}^{L-k+1}c_{j,k}P^{1,\nu}_j(t)\)=\sum\limits_{i=0}^{L-1}\sum\limits_{j=0}^{L-k+1}c_{i,2}c_{j,k}P_i^{1,\nu}(t)P_j^{1,\nu}(t).
$$
By Theorem \ref {positiv}, each polynomial $P_i^{1,\nu}(t)P_j^{1,\nu}(t)$ has non-negative coefficients in terms of $P_k^{1,\nu}$. Then 
$$
\prod\limits_{i=2}^{L}(t-\beta_i)\prod\limits_{i=k}^{L}(t-\beta_i)=\sum\limits_{j=0}^{2L-k}d_jP_j^{1,\nu}(t)
$$
for some $d_j\geq 0$, $j=0,\ldots,2L-k$.
\end {proof}

The polynomial $\Pi_L(t)$ defined in Lemma \ref {beta} is known as the Levenshtein polynomial. We next establish the following two positive definiteness results. 

\begin {lemma}\label {t+1}
Let $n\geq 2$, $k\geq 0$, and $\nu=0,1$. Then the polynomial $(1+t)^\nu P^{1,\nu}_k(t)$ has non-negative coefficients in terms of polynomials $P_i$, $i=0,\ldots,k+\nu$.
\end {lemma}
\begin {proof}
Let $k=0$. If $\nu=0$, we have $(1+t)^\nu P^{1,\nu}_0(t)=1=P_0(t)$ and if $\nu=1$, we have $P_1(t)=t+b$ for some $b\in (-1,1)$ and $(1+t)^\nu P_0^{1,\nu}(t)=1+t=P_1(t)+(1-b)P_0(t)$. That is, $(1+t)^\nu P_0^{1,\nu}(t)$ has non-negative coefficients in terms of $P_i$.

Assume that $k\geq 1$.
If $\nu=0$, then $(1+t)^\nu P_k^{1,\nu}(t)=P_k^{1,0}(t)$ has non-negative coefficients in terms of $P_i^{1,0}$. If $\nu=1$, for every polynomial $Q\in \mathbb P_{k-1}$, 
$$
\int\limits_{-1}^{1}(1+t)P_k^{1,1}(t)Q(t)v_{1,0}(t)\ \! dt=\int\limits_{-1}^{1}P_k^{1,1}(t)Q(t)v_{1,1}(t)\ \! dt=0.
$$
Therefore, $(1+t)P_k^{1,1}(t)=P_{k+1}^{1,0}(t)+a P_k^{1,0}(t)$. Letting $t=-1$, we have $P_{k+1}^{1,0}(-1)+a P_k^{1,0}(-1)=0$. Since all zeros of $P_i^{1,0}$ are simple and located in $(-1,1)$ and its leading coefficient is positive, the sign of $P_i^{1,0}$ is that of $(-1)^i$. Then $a>0$; that is, $(1+t)P_k^{1,1}(t)$ has non-negative coefficients in terms of $P_i^{1,0}$. Since each $P_i^{1,0}$ has non-negative coefficients in terms of $P_j$, see \cite [Equation (3.91)]{Lev}, the assertion of the lemma follows.
\end {proof}

\begin {corollary}\label {P_k}
Let $n\geq 2$, $L\geq 2$, $\nu=0,1$, and numbers $-1<\beta_L<\ldots<\beta_1<1$ be such that the polynomial  
$
\Pi_L(t)
$
is orthogonal to $\mathbb P_{L-2}$ with weight $v_{1,\nu}(t)$. Then the polynomials $(1+t)^\nu$,
$$
(1+t)^\nu \prod\limits_{i=k}^{L}(t-\beta_i)\ \ \ \ \text{and}\ \ \ \ (1+t)^\nu \prod\limits_{i=2}^{L}(t-\beta_i)\prod\limits_{i=k}^{L}(t-\beta_i), \ \ \ 2\leq k\leq L,
$$
have non-negative coefficients in terms of $P_i$.
\end {corollary}
\begin {proof}
If $k=L+1$, the polynomial $(1+t)^\nu$ has non-negative coefficients in terms of $P_i$ by Lemma \ref {t+1}. If $2\leq k\leq L$, the polynomials $(1+t)^\nu \prod\limits_{i=k}^{L}(t-\beta_i)$ and $(1+t)^\nu \prod\limits_{i=2}^{L}(t-\beta_i)\prod\limits_{i=k}^{L}(t-\beta_i)$ have non-negative coefficients in terms of $P_i$ in view of Lemmas \ref {beta} and \ref {t+1}. 
\end {proof}

It remains to prove that the polynomial $(1+t)^\nu \prod\limits_{i=1}^{L}(t-\beta_i)\prod\limits_{i=2}^{L}(t-\beta_i)$ also has non-negative coefficients in terms of~ $P_i$.
To do this, we will use the auxiliary result stated below.
\begin {lemma}\label {n}
Let $n\geq 2$, $\nu\in \{0,1\}$, and $k\geq 0$. Then for every polynomial $p\in \mathbb P_k$, we have
$$
\int\limits_{-1}^{1}(1+t)^\nu p(t)P_k^{1,\nu}(t)v_{0,0}(t)\ \!dt=p(1)\int\limits_{-1}^{1}(1+t)^\nu P_k^{1,\nu}(t)v_{0,0}(t)\ \! dt.
$$
\end {lemma}
\begin {proof}
There exists a polynomial $r\in \mathbb P_{k-1}$ if $k\geq 1$ ($r(t)=0$ if $k=0$) such that $p(t)=(1-t)r(t)+p(1)$. Then
\begin {equation*}
\begin {split}
\int\limits_{-1}^{1}(1+t)^\nu p(t)P_k^{1,\nu}(t)v_{0,0}(t)\ \!dt&=\int\limits_{-1}^{1}(1+t)^\nu (1-t)r(t)P_k^{1,\nu}(t)v_{0,0}(t)\ \!dt+p(1)\int\limits_{-1}^{1}(1+t)^\nu P_k^{1,\nu}(t)v_{0,0}(t)\ \!dt\\
&=\int\limits_{-1}^{1}r(t)P_k^{1,\nu}(t)v_{1,\nu}(t)\ \!dt+p(1)\int\limits_{-1}^{1}(1+t)^\nu P_k^{1,\nu}(t)v_{0,0}(t)\ \!dt\\
&=p(1)\int\limits_{-1}^{1}(1+t)^\nu P_k^{1,\nu}(t)v_{0,0}(t)\ \!dt,
\end {split}
\end {equation*}
which completes the proof.
\end {proof}

Polynomials $P_i$ are related to even degree Gegenbauer polynomials in the following way.  Observe that with $t=2u^2-1$, we have 
$$
v_{\mu,\nu}(t)\ \! dt=2^{\frac {n}{2}+\mu+\nu}u^{2\nu}(1-u^2)^{\mu+\frac {n-3}{2}}\ \! du=\frac {2^{\frac {n}{2}+\mu+\nu}}{\gamma_{n+2\mu}}u^{2\nu}w_{n+2\mu}(u)\ \! du.
$$
\begin {proposition}\label {G'}
Let $n\geq 2$, $k\geq 0$, and $a_{2k}$ be the leading coefficient of the Gegenbauer polynomial $P_{2k}^{(n)}$. Then $P_k(2u^2-1)=\frac {2^k}{a_{2k}}P_{2k}^{(n)}(u)$ and $a_{2k}>0$.
\end {proposition}
\begin {proof}
For every even polynomial $Q\in \mathbb P_{2k-2}$, there is a polynomial $q\in \mathbb P_{k-1}$ such that $Q(u)=q(2u^2-1)$. Then
$$
\int_{-1}^{1}P_k(2u^2-1)Q(u)w_n(u) \ \! du=\frac {\gamma_n}{2^\frac {n}{2}}\int_{-1}^{1}P_k(t)q(t)v_{0,0}(t)\ \! dt
$$
If $Q\in \mathbb P_{2k-1}$ is any odd polynomial, then $P(2u^2-1)$ is orthogonal to it with weight $w_n(u)$, since $w_n(u)$ is even. Every polynomial in $\mathbb P_{2k-1}$ is the sum of an even and an odd polynomial. Therefore, the polynomial $P_k(2u^2-1)$ is orthogonal to $\mathbb P_{2k-1}$ with weight $w_n(u)$. Since $P_k$ is monic, we have $P_k(2u^2-1)=\frac {2^k}{a_{2k}}P_{2k}^{(n)}(u)$. Since $P_{2k}^{(n)}(1)=1>0$ and $P_{2k}^{(n)}$ does not have zeros in $(1,\infty)$, we have $a_{2k}>0$.
\end {proof}

We next establish a positive definiteness result for the polynomial $\Pi_L$.

\begin {lemma}\label {P5}
Let $n\geq 2$, $m=2L+\nu$, where $L\geq 1$ and $\nu\in \{0,1\}$, and numbers $-1<\alpha_1<\ldots<\alpha_m<1$ be symmetric about the origin. Suppose there exists $(C,W)$, a weighted $(m-1,m-1)$-design on $\mathbb{S}^{n-1}$ such that all dot products formed by any pair of points from $C$ lie in the set $A\cup \{-1,1\}$. Denote $\beta_i:=2\alpha_i^2-1$, $i=1,\ldots,L$. Then there exists $\gamma\geq 0$ such that 
$
\Pi_L(t)=\prod\limits_{i=1}^{L}(t-\beta_i)=P_L^{1,\nu}(t)+\gamma P_{L-1}^{1,\nu}(t).
$
\end {lemma}
\begin {proof}
If $L=1$, then there exists $\gamma \in \RR$ such that $\Pi_1(t)=t-\beta_1=P_1^{1,\nu}(t)+\gamma P_0^{1,\nu}(t)$.
In the case $L\geq 2$, let $Q\in \mathbb P_{L-2}$ be an arbitrary polynomial. Then the polynomial $T(u):=\prod\limits_{i=1}^{L}\(u^2-\alpha_i^2\)\cdot (1-u^2)Q(2u^2-1)u^{2\nu}$ has degree at most $4L-2+2\nu=2m-2$ and vanishes on any of the dot products formed by points from $C$. Since $(C,W)$ is a weighted $(m-1,m-1)$-design, by Definition \ref {const} and Remark \ref {rc}, for any fixed point $y_0\in C$, we have
\begin {equation*}
\begin {split}
\int\limits_{-1}^{1}\Pi_L(t)Q(t)v_{1,\nu}(t)\ \! dt &=\frac {2^{\frac {n}{2}+1+\nu+L}}{\gamma_{n+2}}\int\limits_{0}^{1}\prod\limits_{i=1}^{L}\(u^2-\alpha_i^2\)\cdot Q(2u^2-1)u^{2\nu}w_{n+2}(u)\ \! du\\
&=\frac {2^{\frac {n}{2}+\nu+L}}{\gamma_{n}}\int\limits_{-1}^{1}\prod\limits_{i=1}^{L}\(u^2-\alpha_i^2\)\cdot (1-u^2)Q(2u^2-1)u^{2\nu}w_{n}(u)\ \! du=\\
&=\frac {2^{\frac {n}{2}+\nu+L}}{\gamma_{n}}\sum\limits_{y\in C}w_yT(y_0\cdot y)=0.
\end {split}
\end {equation*}
Thus, the degree $L$ polynomial $\Pi_L$ is orthogonal to $\mathbb P_{L-2}$. Since $\Pi_L$ is monic, there exists a constant $\gamma\in \RR$ such that $\Pi_L(t)=P_L^{1,\nu}(t)+\gamma P_{L-1}^{1,\nu}(t)$. 

The rest of the proof is devoted to showing that $\gamma\geq 0$. We assume that $L\geq 1$. If $(C,W)$ is a weighted $(m,m)$-design, then one can take as $Q$ in the above proof any polynomial of degree up to $L-1$. Then the polynomial $T(u)$ will have degree at most $4L+2\nu=2m$ and one will show that $\Pi_L$ is orthogonal to $Q=P_{L-1}^{1,\nu}$ with weight $v_{1,\nu}(t)$; that is, $\gamma=0$.
Therefore, we assume that $(C,W)$ is not a weighted $(m,m)$-design (but it is still an $(m-1,m-1)$-design).

Let $Q\in \mathbb P_m$ be any polynomial. It can be written as
$Q(t)=\sum\limits_{i=0}^{m}p_iP_i(t)$. Then with $P(t):=Q(2t-1)$ taking into account Proposition~\ref {G'}, relation \eqref {equiv1}, and the fact that $a_0=1$, we have
\begin {equation}\label {E_Q}
\begin {split}
E_P(C,W)&=\sum\limits_{i,j=1}^{N}w_iw_jP\((x_i\cdot x_j)^2\)=\sum\limits_{i,j=1}^{N}w_iw_jQ\(2(x_i\cdot x_j)^2-1\)\\
&=\sum\limits_{k=0}^{m}p_k\sum\limits_{i,j=1}^{N}w_iw_jP_k(2(x_i\cdot x_j)^2-1)=\sum\limits_{k=0}^{m}\frac {2^kp_k}{a_{2k}}\sum\limits_{i,j=1}^{N}w_iw_jP_{2k}^{(n)}(x_i\cdot x_j)\\
&=p_0\sum\limits_{i,j=1}^{N}w_iw_j+\frac {2^mp_m}{a_{2m}}\sum\limits_{i,j=1}^{N}w_iw_jP_{2m}^{(n)}(x_i\cdot x_j)=p_0+\frac {2^mp_m}{a_{2m}}\sum\limits_{i,j=1}^{N}w_iw_jP_{2m}^{(n)}(x_i\cdot x_j).
\end {split}
\end {equation}
Denote $b:=\int_{-1}^{1}v_{0,0}(t)\ \! dt$ and let
$$
b_k:=\frac {1}{b}\int\limits_{-1}^{1}(1+t)^\nu P_k^{1,\nu}(t)v_{0,0}(t)\ \! dt,\ \ \ k=0,\ldots,m.
$$
Setting first $Q(t)=(1+t)^\nu \Pi_L(t)=(1+t)^\nu P_{L}^{1,\nu}(t)+\gamma (1+t)^\nu P_{L-1}^{1,\nu}(t)$, we have $p_m=0$, since $Q$ has degree $L+\nu<m$, and from \eqref {E_Q} we have
$$
E_P(C,W)=p_0=\frac {1}{b}\int_{-1}^{1}Q(t)v_{0,0}(t)\ \! dt=b_L+\gamma b_{L-1}.
$$
Since $Q(\beta_i)=0$, $i=1,\ldots,L$, and 
$Q(-1)=0$ when $\nu=1$, we have $Q(2(x_i\cdot x_j)^2-1)=0$ whenever $x_i\cdot x_j\in \{\alpha_1,\ldots,\alpha_m\}$. Then
\begin {equation}\label {theta_m1}
E_P(C,W)=\sum\limits_{i,j=1}^{N}w_iw_jQ\(2(x_i\cdot x_j)^2-1\)=\sum\limits_{i=1}^{N}w_i(w_i+w_i')Q(1)=\theta_{\mu^\ast}Q(1)=2^\nu\theta_{\mu^\ast}\(P_L^{1,\nu}(1)+\gamma P_{L-1}^{1,\nu}(1)\),
\end {equation}
where $\mu^\ast=\sum\limits_{i=1}^{N}w_i\delta_{x_i}$, the quantity $\theta_{\mu^\ast}$ is defined in \eqref {theta_mu}, and $w_i'$ is the weight corresponding to $-x_i$ (if $-x_i\not\in C$ then $w_i'=0$). Then $b_L+\gamma b_{L-1}=2^\nu\theta_{\mu^\ast}\(P_L^{1,\nu}(1)+\gamma P_{L-1}^{1,\nu}(1)\)$; that is,
\begin {equation}\label {2p_L}
b_L-2^\nu\theta_{\mu^\ast}P_L^{1,\nu}(1)+\gamma \(b_{L-1}-2^\nu\theta_{\mu^\ast}P_{L-1}^{1,\nu}(1)\)=0.
\end {equation}
Let now $Q(t)=(1+t)^\nu \Pi_L(t)P_{L}^{1,\nu}(t)$, which has degree $2L+\nu=m$. Then $p_m=1>0$ because $Q$ is monic and, since $(C,W)$ is not a weighted $(m,m)$-design, $\sum\limits_{i,j=1}^{N}w_iw_jP_{2m}^{(n)}(x_i\cdot x_j)>0$. From \eqref {E_Q}, taking into account Lemma \ref {n} and the fact that $a_{2m}>0$ by Proposition~\ref {G'}, we have
$$
E_P(C,W)>p_0=\frac {1}{b}\int\limits_{-1}^{1}(1+t)^\nu \Pi_L(t)P_L^{1,\nu}(t)v_{0,0}(t)\ \!dt=\frac {\Pi_L(1)}{b}\int\limits_{-1}^{1}(1+t)^\nu P_L^{1,\nu}(t)v_{0,0}(t)\ \! dt=\Pi_L(1)b_L.
$$
Since $Q(\beta_i)=0$, $i=1,\ldots,L$, and 
$Q(-1)=0$ when $\nu=1$, similarly to \eqref {theta_m1} we have 
$$
E_{P}(C,W)=\theta_{\mu^\ast}Q(1)=2^\nu\theta_{\mu^\ast}\Pi_L(1)P_L^{1,\nu}(1).
$$
Since $\Pi_L(1)>0$, we have
\begin {equation}\label {b_L}
b_L-2^\nu\theta_{\mu^\ast}P_L^{1,\nu}(1)<0.
\end {equation}
Finally, we let $Q(t)=(1+t)^\nu \(P_{L-1}^{1,\nu}(t)\)^2$, which has degree $2L-2+\nu=m-2$. Then $p_m=0$ and from \eqref {E_Q} taking into account Lemma \ref {n}, we have
\begin {equation*}
\begin {split}
E_p(C,W)&=p_0=\frac {1}{b}\int_{-1}^{1}(1+t)^\nu \(P_{L-1}^{1,\nu}(t)\)^2v_{0,0}(t)\ \! dt\\
&=\frac {P_{L-1}^{1,\nu}(1)}{b}\int_{-1}^{1}(1+t)^\nu P_{L-1}^{1,\nu}(t)v_{0,0}(t)\ \! dt=P_{L-1}^{1,\nu}(1)b_{L-1}.
\end {split}
\end {equation*}
Since $Q(t)\geq 0$, $t\in (-1,1)$, we omit all terms with $x_i\cdot x_j\in \{\alpha_1,\ldots,\alpha_m\}$ in the energy sum below and obtain that
$$
E_P(C,W)=\sum\limits_{i,j=1}^{N}w_iw_jQ\(2(x_i\cdot x_j)^2-1\)\geq \sum\limits_{i=1}^{N}w_i(w_i+w_i')Q(1)=2^\nu\theta_{\mu^\ast}\(P_{L-1}^{1,\nu}(1)\)^2.
$$
The polynomial $P_{L-1}^{1,\nu}$ does not have zeros in $[1,\infty)$. Since its leading coefficient is~$1$, we must have $P_{L-1}^{1,\nu}(1)>0$. Then
$b_{L-1}-2^\nu\theta_{\mu^\ast}P_{L-1}^{1,\nu}(1)\geq 0$. Taking into account \eqref {b_L} from \eqref {2p_L} we obtain that ~$\gamma>0$.
\end {proof}

Lemma \ref {P5} has the following consequence. We agree that $\prod\limits_{i=2}^{1}(t-\beta_i)=1$.

\begin {corollary}\label {2L+1}
Let $n\geq 2$, $m=2L+\nu$, where $L\geq 1$ and $\nu\in \{0,1\}$, and numbers $-1<\alpha_1<\ldots<\alpha_m<1$ be symmetric about the origin. Suppose there exists $(C,W)$, a weighted $(m-1,m-1)$-design on $\mathbb{S}^{n-1}$ such that all dot products formed by any pair of points from $C$ lie in the set $A\cup \{-1,1\}$. Denote $\beta_i:=2\alpha_i^2-1$, $i=1,\ldots,L$.  Then the polynomial  $(1+t)^\nu \prod\limits_{i=2}^{L}(t-\beta_i)\prod\limits_{i=1}^{L}(t-\beta_i)$ has non-negative coefficients in terms of~$P_i$.
\end {corollary}
\begin {proof}
The polynomial $\Pi_L(t)=\prod\limits_{i=1}^{L}(t-\beta_i)$, $L\geq 1$, has non-negative coefficients in terms of $P_i^{1,\nu}$ by Lemma \ref {P5}. In view of Lemma \ref {beta}, the polynomial $\prod\limits_{i=2}^{L}(t-\beta_i)$, $L\geq 2$, also has non-negative coefficients in terms of $P_i^{1,\nu}$ (for $L=1$, this statement is trivial). Then by Theorem~\ref {positiv}, the product $\Pi_L(t)\prod\limits_{i=2}^{L}(t-\beta_i)$ has non-negative coefficients in terms of $P_i^{1,\nu}$. In view of Lemma \ref {t+1}, the polynomial $(1+t)^\nu \Pi_L(t)\prod\limits_{i=2}^{L}(t-\beta_i)$ has non-negative coefficients in terms of $P_i$.
\end {proof}

\begin {theorem}\label {posdefeven}
Let $n\geq 2$, $m\geq 2$, and $A$ be a set of numbers $-1<\alpha_1<\ldots<\alpha_m<1$ symmetric about the origin. Suppose there exists a weighted $(m-1,m-1)$-design $(C,W)$ on $\mathbb{S}^{n-1}$ such that all dot products formed by any pair of points from $C$ lie in the set $A\cup \{-1,1\}$.

Let $g:[-1,1]\to \RR$ be a continuous function such that $g^{(k)}\geq 0$ in $(-1,1)$ for $k=1,\ldots,m-1$, and $p\in \mathbb P_{m-1}$ be the unique polynomial interpolating the values of $g$ at points of the set $B:=\{2\alpha^2_j-1 : j=\OL {1,m}\}$ and derivatives $g'$ at points of $B$ distinct from $-1$. Then the even polynomial $q(t):=p(2t^2-1)$
satisfies
$$
q(t)=\sum\limits_{i=0}^{m-1}d_iP_{2i}^{(n)}(t)
$$
with $d_i\geq 0$, $i=1,\ldots,m-1$. If, in addition, $g^{(k)}>0$ on $(-1,1)$, $k=1,\ldots,m-1$, then $d_i>0$, $i=1,\ldots,m-1$.
\end {theorem}
\begin {proof}
Let $L\geq 1$ and $\nu\in \{0,1\}$ be such that $m=2L+\nu$. For every $t\in [-1,1]$, there exist numbers $\xi_i\in (-1,1)$, $i=\nu,\ldots,m-1$, such that 
\begin {equation}\label {p_i}
p(t)=\nu p(-1)+(t+1)^\nu \sum\limits_{i=\nu}^{L-1+\nu}\frac {g^{(i)}(\xi_i)}{i!}\prod\limits_{j=1}^{i-\nu}(t-\beta_{L+1-j})+(1+t)^\nu\prod\limits_{j=2}^{L}(t-\beta_j)\sum\limits_{i=L+\nu}^{2L-1+\nu}\frac {g^{(i)}(\xi_i)}{i!}\prod\limits_{j=1}^{i+1-L-\nu}(t-\beta_{L+1-j}).
\end {equation}
Since $g^{(k)}\geq 0$ on $(-1,1)$ for each $k=1,\ldots,m-1=2L-1+\nu$, the polynomial $p(t)$ is a linear combination of a constant function, polynomial $(1+t)^\nu$, and polynomials $(1+t)^\nu\prod\limits_{i=k}^{L}(t-\beta_i)$, $k=2,\ldots,L$, and $(1+t)^\nu\prod\limits_{i=2}^{L}(t-\beta_i)\prod\limits_{i=k}^{L}(t-\beta_i)$, $k=1,\ldots,L$, where all polynomials of degree at least $1$ have non-negative coefficients. If $L\geq 2$, by Lemma~\ref {P5}, the polynomial $\Pi_L(t)$ is orthogonal to $\mathbb P_{L-2}$ with weight $v_{1,\nu}(t)$. By Corollaries \ref {P_k} and \ref {2L+1}, polynomials $(1+t)^\nu$,  $(1+t)^\nu\prod\limits_{i=k}^{L}(t-\beta_i)$, $k=2,\ldots,L$, and $(1+t)^\nu\prod\limits_{i=2}^{L}(t-\beta_i)\prod\limits_{i=k}^{L}(t-\beta_i)$, $k=1,\ldots,L$, have non-negative coefficients in terms of $P_i$. If $L=1$, the polynomial $(1+t)^\nu (t-\beta_1)$ has positive coefficients in terms of $P_i$ by Corollary~\ref {2L+1}. The polynomial $(1+t)^\nu$ also has non-negative coefficients in terms of $P_i$, since $P_1$ has a zero in $(-1,1)$. Thus, 
$$
p(t)=d_0+\sum\limits_{i=1}^{m-1}c_iP_i(t),
$$
for some $d_0\in \RR$ and $c_i\geq 0$, $i=1,\ldots,m-1$. By Proposition \ref {G'},
$$
q(u)=p(2u^2-1)=d_0+\sum\limits_{i=1}^{m-1}c_iP_i(2u^2-1)=d_0+\sum\limits_{i=1}^{m-1}d_iP_{2i}^{(n)}(u),
$$
where $d_i:=\frac {2^ic_i}{a_{2i}}\geq 0$, $i=1,\ldots,m-1$.

Assume now that $g^{(k)}>0$ on $(-1,1)$, $k=1,\ldots,m-1$. Then $g^{(i)}(\xi_i)>0$, $i=1,\ldots,m-1=2L-1+\nu$, in~\eqref {p_i}. For each $k=1,\ldots,m-1$, the expansion in terms of $P_i$ of the product of factors $(t-\beta_j)$ and $(1+t)^\nu$ of degree $k$ in \eqref {p_i} contains $P_k$ with coefficient $1$. The expansion of every such product of degree greater than $k$ contains $P_k$ with a non-negative coefficient. Then $c_k$ is strictly positive as the sum of these coefficients multiplied by strictly positive values~$\frac {g^{(i)}(\xi_i)}{i!}$. Therefore, $d_i>0$, $i=1,\ldots,m-1$.
\end {proof}


\noindent{\bf Acknowledgment.}  The research of the second author was supported, in part, by Bulgarian NSF grant KP-06-N72/6-2023. 
The research of the third author is supported, in part, by the Lilly Endowment.
The research of the sixth author was supported, in part, by the European Union-NextGenerationEU, through the National Recovery and Resilience Plan of the Republic of Bulgaria,
Grant no. BG-RRP-2.004-0008.


\begin{thebibliography}{0}
\bibitem{A2} Andreev, N. N., A minimal design of order $11$ on the $3$-sphere,
Math. Notes {\bf 67}, 417--424 (2000).





\bibitem {BilGlaPar2021}
Bilyk, D., Glazyrin, A., Matzke, R., Park, J., Vlasiuk, O., Energy on spheres and discreteness of minimizing measures, Journal of Functional Analysis {\bf 280} (2021), no. 11: 108995.

\bibitem {BilGlaPar2022}
Bilyk, D., Glazyrin, A., Matzke, R., Park, J., Vlasiuk, O., Optimal measures for $p$-frame energies on spheres, Rev. Mat. Iberoam. {\bf 38} (2022), no. 4, pp. 1129--1160.




\bibitem{B2} Borodachov, S. V., Polarization problem on a higher-dimensional sphere for a simplex, Discr. Comput. Geom. {\bf 67}, 525--542 (2022). 

\bibitem{BorESI}
S.V. Borodachov, Min-max polarization for certain classes of sharp configurations on the sphere, {\it Workshop "Optimal Point Configurations on Manifolds"}, ESI, Vienna, January 17--21, 2022. https://www.youtube.com/watch?v=L-szPTFMsX8



\bibitem{B} Borodachov, S. V., Min-max polarization for certain classes of sharp configurations on the sphere, Constr. Approx. (2023), 
https://doi.org/10.1007/s00365-023-09661-1.

\bibitem{Bor-new} Borodachov, S. V., Absolute minima of potentials of certain regular spherical configurations, 
J. Approx. Theory, {\bf 294}, art. 105930 (2023). 

\bibitem{Bor-2} Borodachov, S. V., Absolute Minima of Potentials of a Certain Class of Spherical Designs, arXiv:2212.04594.

\bibitem{Bor-FL} Borodachov, S. V., Odd strength spherical designs attaining the Fazekas–Levenshtein bound for covering and universal minima of potentials,
Aequat. Math., {\bf 98}, 509--533 (2024).
\bibitem {BorBoyDraenergy}
Borodachov, S. V., Boyvalenkov, P. G., Dragnev, P. D., Hardin, D. P., Saff, E. B., Stoyanova, M. M., Energy bounds for weighted spherical codes and designs via linear programming, arXiv: 2403.07457.

\bibitem{BHS} Borodachov, S. V., Hardin, D. P., Saff, E. B., Discrete Energy on Rectifiable Sets, Springer Monographs in Mathematics, Springer, 2019.




\bibitem{BDHSS} Boyvalenkov, P., Dragnev, P., Hardin, D., Saff, E., Stoyanova, M., Universal lower bounds for potential energy of spherical codes,
Constr. Approx. {\bf 44}, 385--415 (2016).

\bibitem{BDHSS-AMP} Boyvalenkov, P., Dragnev, P., Hardin, D., Saff, E., Stoyanova, M., 
Energy bounds for codes in polynomial metric spaces, Anal. Math. Phys. {\bf 9}, 781--808 (2019).

\bibitem{BDHSS-DCC19} Boyvalenkov, P., Dragnev, P., Hardin, D., Saff, E., Stoyanova, M., 
On spherical codes with inner products in a prescribed interval, Des. Codes Crypt. {\bf 87}, 299--315 (2019).

\bibitem{BDHSSMathComp} Boyvalenkov, P., Dragnev, P., Hardin, D., Saff, E., Stoyanova, M., 
Bounds for spherical codes: the Levenshtein framework lifted, Math. Comp. {\bf 90}, 1323--1356 (2021).

\bibitem{BDHSS-JMAA} Boyvalenkov, P., Dragnev, P., Hardin, D., Saff, E., Stoyanova, M.,  On polarization of spherical codes and designs, JMAA, {\bf 524},  art. 127065 (29 pages), 2023. 
.

\bibitem{BDHSS-Sharp} Boyvalenkov, P., Dragnev, P., Hardin, D., Saff, E., Stoyanova, M.,  Universal minima of discrete potentials for sharp spherical codes, submitted (arXiv:2211.00092).







\bibitem{CohKum2007} Cohn, H., Kumar, A., Universally optimal distribution of points on spheres,
J. Amer. Math. Soc. {\bf 20}, 99--148 (2007).




\bibitem{DGS} Delsarte, P., Goethals, J.-M., Seidel, J. J., Spherical codes and designs, 
Geom. Dedic. {\bf 6}, 363--388 (1977).






\bibitem{G} Gasper, G., Linearization of the product of Jacobi polynomials I, Canadian Journal of Mathematics, {\bf 22}, no. 1, 171--175 (1970).
 

\bibitem {GoeSei1981}
Goethals, J. M., Cubature formulas, polytopes, and spherical designs. In: Davis, C., Gr\"unbaum, B., Sherk, F.A. (eds) The Geometric Vein. Springer, New York, NY, 203--218.


\bibitem {HarKenSaf2013}
Hardin, D., Kendall, A., Saff, E.,
Polarization optimality of equally spaced points on the circle for discrete potentials.
{\it Discrete Comput. Geom.} {\bf 50} (2013), no. 1, 236--243. 









\bibitem{Lev} Levenshtein, V. I., Universal bounds for codes and designs, in Handbook of Coding Theory, V. S. Pless and W. C. Huffman, Eds., 
Elsevier, Amsterdam, Ch.~6, 499--648 (1998).





\bibitem{NR1} Nikolov, N., Rafailov, R., \emph{On the sum of powered distances to certain sets of points on the circle}, Pacific J. Math. {\bf 253}, no. 1, 157--168 (2011).

\bibitem{NR2} Nikolov, N., Rafailov, R., On extremums of sums of powered distances to a finite set of points, Geom. Dedic. {\bf 167}, 69--89 (2013).






\bibitem {Sto1975circle}
Stolarsky, K., The sum of the distances to certain pointsets on the unit circle, {\it Pacific J. Math.} {\bf 59} (1975), no. 1, 241--251. 
\bibitem {Sto1975}
Stolarsky, K., The sum of the distances to $N$ points on a sphere, {\it Pacific J. Math.} {\bf 57} (1975), no. 2, 563--573.




\bibitem{Sze1975} Szeg\"{o}, G., Orthogonal polynomials, AMS Col. Publ., {\bf 23}, Providence, RI, 1939.





\bibitem{Y} Yudin, V. A.,
Minimal potential energy of a point system of charges,
Discr. Math. Appl., {\bf 3} (1983), 75--81.

\bibitem {Yud2005}
Yudin, V. A.,
Distribution of the points of a design on a sphere
Izv. Math., {\bf 69} (2005), no. 5, 1061--1079


\end{thebibliography}
\end{document}